\documentclass[12pt]{article}
\title{A holomorphic point of view\\ about geodesic completeness
}
\author{Claudi Meneghin
}

\begin{document}
\maketitle
\bibliographystyle{plain} 
\parindent=8pt

\begin{abstract}
\noindent
We propose
to apply the idea of
analytical continuation
in the complex domain to
the problem of
geodesic
completeness.
We shall analyse rather in detail
the cases of analytical warped
products of real lines, these
ones in parallel with their
complex
counterparts,
and of Clifton-Pohl torus,
to show that our definition
sheds a bit of new light on
the behaviour of 'singularities'
of geodesics in space-time.
We also show that some
geodesics, which 'end' at
finite time in the classical
sense, can be naturally
continued besides their ends.
As a matter of fact, complex
metrics naturally show
a meromorphic behaviour,
or a degenerating one, so
we shall study also this fact
in detail.
\end{abstract}

%\documentclass[10pt]{article}
%\begin{document}
\font\sdopp=msbm10
\def\CI {\sdopp {\hbox{C}}}

\def\quan{\vrule height6pt width6pt depth0pt}
\def\QUAN{\nobreak $\ $\quan\hfill\vskip0.6truecm\par}
\def\BETA{\mathop{\beta}\limits}
\def\GAMMA{\mathop{\gamma}\limits}
\def\VI{\mathop{v}\limits}
\def\UI{\mathop{u}\limits}
\def\VII{\mathop{V}\limits}
\def\WI{\mathop{w}\limits}
\def\ZETA{\mathop{Z}\limits}

\newtheorem{definition}{Definition}%[chapter]
\newtheorem{lemma}[definition]{Lemma}
\newtheorem{proposition}[definition]{Proposition}
\newtheorem{theorem}[definition]{Theorem}
\newtheorem{corollary}[definition]{Corollary}
\newtheorem{remark}[definition]{Remark}

\font\sdopp=msbm10
\def\ESSE {\sdopp {\hbox{S}}}
\def\ERRE {\sdopp {\hbox{R}}}
\def\CI {\sdopp {\hbox{C}}}
\def\ENNE{\sdopp {\hbox{N}}}
\def\ZETA{\sdopp {\hbox{Z}}}
\def\PI {\sdopp {\hbox{P}}}
\def\P{\hbox{\boldmath{}$P$\unboldmath}}
\def\Y{\hbox{\boldmath{}$Y$\unboldmath}}
\def\tr{\hbox{\boldmath{}$tr$\unboldmath}}
\def\f{\hbox{\boldmath{}$f$\unboldmath}}
\def\u{\hbox{\boldmath{}$u$\unboldmath}}
\def\v{\hbox{\boldmath{}$v$\unboldmath}}
\def\U{\hbox{\boldmath{}$U$\unboldmath}}
\def\V{\hbox{\boldmath{}$V$\unboldmath}}
\def\W{\hbox{\boldmath{}$W$\unboldmath}}
\def\id{\hbox{\boldmath{}$id$\unboldmath}}
\def\alph{\hbox{\boldmath{}$\aleph$\unboldmath}}
\def\bet{\hbox{\boldmath{}$\beta$\unboldmath}}
\def\gam{\hbox{\boldmath{}$\gamma$\unboldmath}}
\def\U{\mathop{u}\limits}
\def\f{\hbox{\boldmath{}$f$\unboldmath}}
\def\g{\hbox{\boldmath{}$g$\unboldmath}}
\def\h{\hbox{\boldmath{}$h$\unboldmath}}

\def\TTT{}
\def\obmet#1{}
\sloppy

\bibliographystyle{plain}
\parindent=8pt
%caratteri cirillici
%\font\cir=wncyb10
%\def\Iu{\cir\hbox{YU}}
%\def\Ze{\cir\hbox{Z}}
%\def\pe{\cir\hbox{P}}
%\def\Ef{\cir\hbox{F}}
%
%\font\sdopp=msbm10
%\def\DI{\sdopp{\hbox{D}}}
%\def\ESSE {\sdopp {\hbox{S}}}
%\def\ERRE {\sdopp {\hbox{R}}}
%\def\CI {\sdopp {\hbox{C}}}
%\def\ENNE{\sdopp {\hbox{N}}}
%\def\ZETA{\sdopp {\hbox{Z}}}
%\def\PI {\sdopp {\hbox{P}}}
%\def\M{\hbox{\boldmath{}$M$\unboldmath}}
%\def\N{\hbox{\boldmath{}$N$\unboldmath}}
%\def\P{\hbox{\boldmath{}$P$\unboldmath}}
%\def\Y{\hbox{\boldmath{}$Y$\unboldmath}}
%\def\tr{\hbox{\boldmath{}$tr$\unboldmath}}
%\def\f{\hbox{\boldmath{}$f$\unboldmath}}
%\def\u{\hbox{\boldmath{}$u$\unboldmath}}
%\def\v{\hbox{\boldmath{}$v$\unboldmath}}
%\def\U{\hbox{\boldmath{}$U$\unboldmath}}
%\def\V{\hbox{\boldmath{}$V$\unboldmath}}
%\def\W{\hbox{\boldmath{}$W$\unboldmath}}
%\def\id{\hbox{\boldmath{}$id$\unboldmath}}
%\def\alph{\hbox{\boldmath{}$\aleph$\unboldmath}}
%\def\bet{\hbox{\boldmath{}$\beta$\unboldmath}}
%\def\gam{\hbox{\boldmath{}$\gamma$\unboldmath}}
%\def\U{\mathop{u}\limits}
%\def\f{\hbox{\boldmath{}$f$\unboldmath}}
%\def\g{\hbox{\boldmath{}$g$\unboldmath}}
%\def\h{\hbox{\boldmath{}$h$\unboldmath}}
\def\IM{\hbox{\boldmath{}$i$\unboldmath}}

\def\Ch{\hbox{\rm Ch}}
\def\CIRC{\mathop{\tt o}\limits}

\def\oi{\mathop{\omega}\limits}
\def\ei{\mathop{\eta}\limits}
\def\FI{\mathop{\varphi}\limits}

\def\BBB{\sl}
%\def\quan{\vrule height6pt width6pt depth0pt}
%\def\QUAN{\nobreak $\ $\quan\hfill\vskip0.6truecm\par}
%\def\BETA{\mathop{\beta}\limits}
%\def\GAMMA{\mathop{\gamma}\limits}
%\def\VI{\mathop{v}\limits}
%\def\UI{\mathop{u}\limits}
%\def\VII{\mathop{V}\limits}
%\def\WI{\mathop{w}\limits}
%\def\ZETA{\mathop{Z}\limits}

%font\sdopp=msbm10
%def\ESSE {\sdopp {\hbox{S}}}
%def\ERRE {\sdopp {\hbox{R}}}
%def\CI {\sdopp {\hbox{C}}}
%def\ENNE{\sdopp {\hbox{N}}}
%def\ZETA{\sdopp {\hbox{Z}}}
%def\PI {\sdopp {\hbox{P}}}
%def\M{\hbox{\boldmath{}$M$\unboldmath}}
%def\N{\hbox{\boldmath{}$N$\unboldmath}}
%def\P{\hbox{\boldmath{}$P$\unboldmath}}
%def\Y{\hbox{\boldmath{}$Y$\unboldmath}}
%def\tr{\hbox{\boldmath{}$tr$\unboldmath}}
%def\f{\hbox{\boldmath{}$f$\unboldmath}}
%def\u{\hbox{\boldmath{}$u$\unboldmath}}
%def\v{\hbox{\boldmath{}$v$\unboldmath}}
%\def\U{\hbox{\boldmath{}$U$\unboldmath}}
%\def\V{\hbox{\boldmath{}$V$\unboldmath}}
%\def\W{\hbox{\boldmath{}$W$\unboldmath}}
%\def\id{\hbox{\boldmath{}$id$\unboldmath}}
%\def\alph{\hbox{\boldmath{}$\aleph$\unboldmath}}
%\def\bet{\hbox{\boldmath{}$\beta$\unboldmath}}
%\def\gam{\hbox{\boldmath{}$\gamma$\unboldmath}}
%\def\U{\mathop{u}\limits}
%\def\f{\hbox{\boldmath{}$f$\unboldmath}}
%\def\g{\hbox{\boldmath{}$g$\unboldmath}}
%\def\h{\hbox{\boldmath{}$h$\unboldmath}}

\font\sdopp=msbm10
\def\DI {\sdopp {\hbox{D}}}

\def\M{\hbox{\tt\large M}}
\def\N{\hbox{\tt\large N}}
\def\T{\hbox{\tt\large T}}

\def\e{\hbox{\boldmath{}$e$\unboldmath}}

\def\labelle #1{\label{#1}}

\section{Foreword}
Within
the framework of
Riemannian geometry,
geodesic
and metric completeness are
well known to be equivalent:
this is Hopf-Rinow's theorem, a consequence of
 the positivity
of Riemannian metrics.

This equivalence does not hold
for semi-Riemannian metrics,
and there exist even compact
lorentzian manifolds which are
geodesically incomplete: a
 well known exemple is
Clifton-Pohl torus (see e.g.
\cite{oneill}, 7.16).

In this paper
we propose a
definition of
geodesic completeness
from a complex point of view, that is to say
we shall look rather at complexified
pseudo-Riemannian manifolds with
complex-symmetric metrics.

By a philosophical point of
view,
our goal is to shed a little
bit of  light on the
behaviour of 'singularities'
of geodesics in space-time and
show that some geodesics,
which seem to 'end' at finite
time can be naturally
continued besides their ends.
 This will be done by
running along complex trips
close
to the real line.

Since our approach will use
complex-analytical methods
and analytical continuation
leads in
general to poles and zeroes,
we shall need the idea
of a {\it meromorphic metric}
on a complex manifold $\M$
(see \cite{lebrun} for
the definition of a
{\it holomorphic metric};
see also \cite{manin};
see \cite{bfv} about
the relationship with
anti-K\"ahler geometry; see
\cite{cassab} and \cite{cassac} for physical motivations).
This will amount to a
possibly degenerating symmetric
meromorphic section
of
the twice covariant holomorphic tensor bundle
 ${\cal T}_{0}^{2}\M$.

Of course, it carries
no 'signature'.
However, by simmetry,
it induces
a canonical meromorphic
Levi-Civita's connexion
on $\M$, allowing
 to define geodesics as the
auto-parallel paths.
For the sake of completeness
this aspects will be
dealt with in some details.

It is worth noticing that,
if $\M$ arises as a
'complexification'
of a semi-Riemannian
manifold $\N$,
it is easily seen that
the real geodesics  of $\N$
are restrictions to
the real axis of the complex
ones of $\M$
and vice versa
(see \cite{lebrun}).

This fact sometimes allows us to 'flank' isolated singularities on the real line, i.e. to 'connect'
geodesics which, in the usual sense,
are completely unrelated.

The problem of lorentzian geodesic completeness
is investigated in \cite {beemerlich}, \cite{oneill},
\cite{carri},
\cite{romero}.

We suggest a more formal idea
of
our notion of {\sl completeness}
(see also definition \ref{completessa})
:
given a complexification $d:\N\rightarrow\M$
and a real analytic curve $\gamma:[a,b]\rightarrow\N$,
$\gamma$ will be
told to be
{\sl complete}
provided that
$d\circ\gamma$ can be continued to all points
in the real line,  with at most
a discrete set of exceptional values,
taking 'real values' (i.e. in $d(\N)$).

Finally,
we report the main existence-and-uniqueness
theorem of ordinary differential equation theory in the complex domain:
let $W_0
$ be a complex $N-$tuple,
$z_0\in\CI$;
let $F$ be a $\CI^N-$valued holomorphic mapping in
$\prod_{j=1}^N\DI
(W_0^j,b
)\times\DI
(z_0,a
)$,
($a,b\in\ERRE$)
with $C^0-$norm $M$ and
$C^0-$norm of each
${\partial F}/{\partial w^j}$ ($j=1..N$) not exceeding $K\in\ERRE$.

\begin{theorem}
If
$r<min(a,b/M,1/K)$,
there exists a unique
holomorphic mapping
$
W\colon \DI
(z_0,r
)$ $\rightarrow
\prod_{j=1}^N\DI
(W_0^j,b
)
$
such that
$W^{\prime}=F(W(z),z)$ and
$W(z_0)=W_0$.
(see e.g. \cite{hille}, th 2.2.2, \cite{ince} p.281-284)
\label{gerexuq}
\end{theorem}

NB: in the following we shall abbreviate 'holomorphic function element',
resp. 'holomorphic function germ' by HFE resp. HFG.

\section{Analytical
continuation}
In the following,  ${\cal U}$ will be a region in the complex plane and $\M$
a complex manifold: the idea of the analytical continuation of a holomorphic mapping element
(or of a germ)
$f:{\cal U}\rightarrow\M$  is well known
see e.g. \cite{cassa}, chap. 5, rather than \cite{forster}, 1.7-1.8)
and amounts to a quadruple $Q_{\M}=(S,\pi,j,F)$, where
$S$ is a connected Riemann surface over a region of $\CI$,
$\pi\,\colon\, S\rightarrow \CI$ is a nonconstant holomorphic mapping
such that $U\subset \pi(S)$,
$j\,\colon\, U\rightarrow S$ is a holomorphic  immersion such that $\pi\circ j=id\vert_{U}$
and
$F\,\colon\, S\rightarrow \M$ is a holomorphic mapping such that $F\circ j=f$.
Each finite branch point is
kept into account by the fact of lying 'under' some critical point of $\pi$: for example, the Riemann surface of $\sqrt{z}\vert_{D(1,1/2)}$, with $\sqrt{1}=1$,
is $\left(\CI, \zeta\mapsto\zeta^2,
\sqrt{z}\vert_{D(1,1/2)},
\zeta\mapsto\zeta
\right)$: the double branch point $z=0$ of the analytical continuation lies under the critical point $\zeta=0$, branching being taken into account by the squaring function.

A {\sf morphism} between two analytical continuations
$
\left(S,\pi,j,F\right)
$
and
$\left(T,\varrho,\ell,G\right)
$
of the same element $\left(U,f    \right)$
is a holomorphic mapping $h\,\colon\, T\rightarrow S$
such that $h\circ\ell=j$.
Note that
a morphism between two analytical
continuations is a nonconstant (in particular
open) mapping, uniquely determined on $j(U)$, hence everywhere on $S$, by $\ell\circ j^{-1}$.
Moreover, $\varrho\circ h=\pi$ and
$G\circ h=F$ on $j(U)$, hence everywhere on $S$.

The only existing morphism between an
analytical continuation and itself is the identity mapping; the composition of two morphisms is
a morphism; if a morphism admits a
holomorphic inverse mapping, this is
 a morphism too: in such a case we talk about {\sf isomorphisms} of analytical continuations.
%%%%%%%%%%%%%%%%%%%%

\begin{definition}
An analytical continuation
$S$
 of the element $\left(U,f    \right)$
is {\sf  maximal} if for every regular analytical continuation $\widehat S$
 of $\left(U,f    \right)$ there exists a morphism $h\,\colon\, S\rightarrow \widehat S$.
\end{definition}
Two maximal continuations of the same element must be isomorphic, since they admit morphisms one into each other; thus the maximal regular  analytical continuation of an
 element is unique up to isomorphisms.

The following is a well known
result (we refer to \cite{cassa}, th.5.6.4, pages 262-266).

\begin{theorem}
Every element $\left(U,f    \right)$ (hence every germ) of holomorphic function
admits a maximal analytical continuation,
called the {\BBB Riemann surface}, of $\left({\cal U},f   \right)$.
\labelle{ex_max_sup}
\end{theorem}
\begin{definition}
A
{\BBB  logarithmic singularity} $q$ of
$Q_{\M}=\left(S,\pi,j,F\right)
$  (in the following: {\BBB L-singularity}) is a decreasing sequence of
open sets
$\{V_k\}_{k\geq K}$ of $S$
such that:\\
$\bullet$ (LS1)
for every $k\geq K$,
 $V_k$ is a connected component of
${\pi^{-
1}(D(z_0,\frac{1}{k})\setminus \{z_0\})}$
and $\pi\vert_{V_k}$ is a topological
covering of $(D(z_0,\frac{1}{k})$;\\
$\bullet$ (LS2)
${\bigcap_{k\geq K}\overline{V_k}=
\emptyset}$;\\
$\bullet$ (LS3)
for every $k\geq K$ and every
(real) nonconstant
closed path $\gamma:[0,1]
\rightarrow D(z_0,1/k)
\setminus\{z_0\}$, with
nonzero winding number
around $z_0$, every
lifted path $\tilde\gamma:[0,1]
\rightarrow \pi^{-1}
(D(z_0,1/k)\setminus\{z_0\})$
 is such that
$\tilde\gamma(0)\not=\tilde\gamma(1)$;\\
$\bullet$ (LS4)
there exists $m\in\M$ such that
${\bigcap_{k\geq K}\overline{F(V_k)}=
m}$.
%
%\label{log_singol2}
\label{log_sing}
%\label{sps}
\end{definition}

Consider now the set $B$ of the
 L-singularities of $Q_{\M}$: let $S^{\sharp}:=S\bigcup B$ as a set and
introduce a topology on $S^{\sharp}$:
open sets are the open sets in $S$ and a fundamental neighbourhood system of
the L-singularity $q={\{V_k\}_{k\geq K}}\in B$ is yielded by the sets $V_k^{\sharp}=V_k\bigcup \{q\}$.
It is easily seen that
$S^{\sharp}$ admits no complex structure at
 $q={\{V_k\}_{k\geq K}}$.
Indeed, were there one, we could find charts $({\cal W},\phi   )$ around $q$ and $({\cal V},\psi   )$ around $z_0$ such that
$
\psi\circ\pi\circ\phi^{-1}(\zeta   )
=\zeta^N
$
for some integer $N>0$.
This fact would imply $\pi\vert_{{\cal W}\setminus\{q\}}$ to be a n-sheeted covering of ${\cal V}\setminus\{z_0\}$;
it is easily seen tha this fact would contradict
(LS2) in definition \ref{log_sing}.

\begin{lemma}
{\bf (A)}: $\pi$ admits a unique continuous extension $\pi^{\sharp}$ to $S^{\sharp}$;
{\bf (B)}:  for every logarithmic
singularity $r$ of $Q_{\M}$, $F$ admits a unique continuous extension
 $F^{\sharp}$ to $r$.
\label{funest}
\label{proest}
\end{lemma}
{\bf Proof:} {\bf (A)}: let $b\in B$ and $\{V_k\}$ be the sequence spotting $b$: define
$
\pi^{\sharp}(q)=\pi(q)$ if $q\in V_k$
and
$
\pi^{\sharp}(b)=z_0
$,
where $z_0$ is the common centre of the discs onto which the $V_k's$ are projected.
Now $\pi^{\sharp}$ is continuous at all points in $V_k$; moreover, for every neighbourhood $G$ of $z_0$,
$\pi^{\sharp\ -1}(G)\supset \pi^{\sharp\ -1}(z_0) \bigcup\pi^{-1}(G\setminus\{z_0\})$,
hence, if we set  $H=\{b\}\bigcup\pi^{-1}(G\setminus\{z_0\})$, we have that $H$ is a neighbourhood of $b$ in $S^{\sharp}$ such that $\pi^{\sharp}(H)\subset G$, proving continuity at $b$. Arguing by density, we conclude that this extension is unique; the proof of {\bf (B)} is analogous.
\QUAN
\begin{definition}
A quadruple $ Q_{\M}^{\natural}
=(S^{\natural},\pi^{\natural},j^{\natural},
F^{\natural})$, is
an {\BBB analytical continuation with L-singularities} of the function element $(U,f)$ if there exists an analytical continuation
$ Q_{\M}$ of
$(U,f)$ such that
$ S^{\natural}\setminus S$ consists of L-singularities of $F$, $\pi^{\natural}$ is the unique continuous extension of $\pi$ to
 $S^{\natural}$,
$j^{\natural}=id_{S\rightarrow S^{\natural}}\circ j$ and
$F$ admits a unique continuous
extension $ F^{\natural}$ to
$S^{\natural}$.
$Q_{\M}^{\natural}$ is:
{\BBB maximal} provided that
so is $Q_{\M}$ and
$Q_{\M}^{\natural}\setminus
 Q_{\M}$ contains all
L-singularities of $Q_{\M}$.
\label{loganalcont}
\label{ext_rsb}
\end{definition}

\begin{lemma}
\label{inverse}
{\bf (1)}: let $\f$ and $\g$ be two
complex-valued
holomorphic germs,
each one inverse of the other;
 let $(   R,\pi,j,F)$ and $(S,\rho,\ell,G)$
be their respective Riemann surfaces: then $F(R)=\rho(S)$;\\
{\bf (2)}: \label{quasiinverse}
let
$\f$, $\g$, $\h$
be three HFG's such that
$\f\circ \g=\h$.
Let $(   R,\pi,j,F)$
be the Riemann surface of
$\f$, $(S,\rho,\ell,G)$
the one of $\g$ and
$(T,\sigma,m,H)$
the Riemann surface with L-singularities of $\h$: then $F(R)\setminus(\CI\setminus(
\sigma(T)
   )   )\subset\rho(S)$.
\end{lemma}
{\bf Proof:}
{\bf (1)}
{\tt\large a):} $F(R)\subset\rho(S)$: let $\xi\in R$ and $F(\xi)=\eta$; there exist:
an open neighbourhood ${\cal U}_1$ of $\xi$;
open subsets ${\cal U}_2\subset\pi({\cal U}_1   )$ and ${\cal V}_2\subset F({\cal U}_1   )$ and
a biholomorphic function $g_2:{\cal V}_2\rightarrow{\cal U}_2$, with inverse function $f_2:{\cal U}_2\rightarrow{\cal V}_2$
{such that}:
$({\cal U}_2,f_2   )$ and $({\cal U},f  )$ are connectible and so are
$({\cal V}_2,g_2   )$ and $({\cal V},g  )$.
By construction there hence exist two holomorphic immersions
$\widetilde{\hbox{\j}} :{\cal U}_2\rightarrow R\hbox{ and }
\widetilde\ell:{\cal V}_2\rightarrow S$
such that $\pi\circ\widetilde{\hbox{\j}}=\id$ and $\rho\circ\widetilde\ell=\id$.
Let ${\cal V}_1=F(U)_1$ and
$
\Sigma=\{(x,y   )\in{\cal U}_1\times{\cal V}_2: F(x)=y\}
$;
moreover let $J:{\cal V}_2\rightarrow\Sigma$
be defined by setting $ J(v)=(\widetilde{\hbox{\j}}\circ  g_2(v),v )$.
Then $(\Sigma,pr_2,J,\pi\circ pr_1   )$ is an analytical continuation of $({\cal V}_2,g_2   )$; indeed $ \pi\circ pr_1\circ J=\pi\circ\widetilde{\hbox{\j}}\circ g_2=g_2$. But $({\cal V_2},g_2   )$ is connectible with $({\cal V},g     )$, hence
$ (\Sigma,pr_2,J,\pi\circ pr_1   )$ is an analytical continuation of $({\cal V},g )$.
Eventually, there  exists a holomorphic function $h:\Sigma\rightarrow S$ such that $\rho\circ h=pr_2$: hence
$
\eta=pr_2(\xi,\eta   )=\rho\circ h
(\xi,\eta   )\in\rho(S   )
$.\\
{\tt\large b):} $\rho(S)\subset F(R)$: let $s\in S$: there is a neighbourhood  $V$ of $s$ in $S$ such that $V\setminus\{s\}$ consists entirely of regular points both of $\rho$ and $G$, not excluding that $s$ itself be regular for $\rho$ or $G$ or both.
This fact means that for each $s^{\prime}\in V\setminus\{s\}$ there exists a HFE $({\rho(s^{\prime})},{\cal V}^{\prime},\widetilde g_{s^{\prime}}   )$ connectible with $({\cal V},g   )$ and, besides, a holomorphic immersion $\widetilde\ell:{\cal V}^{\prime}\rightarrow V$.
By {\tt\large a):} already proved, $G(s)\in\pi(R)$, hence there exist
$p\in R$ such that $\pi(p)=G(s)$ and
a neighbourhood $W$ of $p$ in $R$ such that $ \pi^{-1}(\widetilde g({\cal V}^{\prime}   )   )\bigcap W\not=\emptyset$.
Set
$ W^{\prime}=\pi^{-1}(\widetilde g({\cal V}^{\prime}   )   )\bigcap W$:
we may suppose, without loss of generality, that $\pi$ is invertible on $W^{\prime}$: hence there exists a (open) holomorphic immersion $\widetilde{\hbox{\j}}:\widetilde g({\cal V}^{\prime}   )\rightarrow W$.
Therefore, for each $\zeta\in\widetilde{\hbox{\j}}(\widetilde g({\cal V}^{\prime}   )   )$, there exists $ \eta\in\widetilde\ell({\cal V}^{\prime}   )$ such that $ F(\zeta)=F(\widetilde{\hbox{\j}} \circ \widetilde g\circ \rho(\eta)  )$.
Now, by definition of analytical continuation there holds $ F\circ\widetilde{\hbox{\j}}\circ\widetilde g=\id$, hence we have $ F(\zeta)=\rho(\eta)$.
Consider now the holomorphic function
$\Xi:W\times V\rightarrow\CI$
defined by setting $ \Xi(w,v   )=F(w)-\rho(v)$: we have
$
\Xi
\equiv 0
$ on the open set
$ {\widetilde{\hbox{\j}}(\widetilde g({\cal V}^{\prime}   )   )\times\widetilde\ell({\cal V}^{\prime})}$,
thus $\Xi\equiv 0$ on $W\times V$:
this in turn implies $ F(p)=\rho(s)$.
Therefore we have proved that for each $s\in S$ there exists $p\in R$ such that $ F(p)=\rho(s)$: this eventually implies that $\rho(S)\subset F(R)$.\\
{\bf (2)}:
 let $\xi\in R$ such that
 $\displaystyle\eta=F(\xi)
\not\in\CI
\setminus\left(
\sigma(T)
   \right)  $:
there exist:
an open neighbourhood
 ${\cal U}_1$ of $\xi$,
open subsets ${\cal U}_2
\subset\pi\left({\cal U}_1
  \right)$, ${\cal V}_2
\subset F\left({\cal U}_1
\right)$ and  ${\cal W}_2
\subset \sigma\left(T  \right)$
and
biholomorphic functions
$\displaystyle f_2:{\cal U}_2\rightarrow{\cal W}_2$,
$\displaystyle g_2:{\cal V}_2\rightarrow{\cal U}_2$ and
$\displaystyle h_2:{\cal V}_2\rightarrow{\cal W}_2$
{\TTT such that}:
$\left({\cal U}_2,f_2
  \right)$ and $\f$
are connectible,
$\left({\cal V}_2,g_2
  \right)$ and $\g$ are
 connectible,
$\left({\cal V}_2,h_2
 \right)$ and $\h$
are connectible,
and $f_2\circ g_2=h_2$.

By construction there hence
exist three holomorphic
immersions
$
\widetilde{\hbox{\j}}:{\cal U}_
2\rightarrow R$,
$
\widetilde\ell:{\cal V}_2
\rightarrow S$,
$
\widetilde m:{\cal W}_2
\rightarrow T
$
such that
$\displaystyle\pi\circ
\widetilde{\hbox{\j}}=\id$,
$\displaystyle\rho\circ
\widetilde\ell=\id$ and
$\displaystyle\sigma\circ
\widetilde m=\id$.

Let ${\cal V}_1=F\left(U_1\right)$, ${\cal W}_1$ be the connected component of $\displaystyle \sigma^{-1}\left(F\left({\cal U}_1   \right)   \right)$ in $T$ and
$
\Sigma=\{\left(x,y   \right)\in{\cal U}_1\times{\cal W}_1: F(x)=H(y)\}
$;
moreover let $J:{\cal V}_2\rightarrow\Sigma$
be defined by setting $\displaystyle J(v)=\left(\widetilde{\hbox{\j}}\circ  g_2(v),\widetilde m(v) \right)$.

Then $\displaystyle \left(\Sigma,pr_2,J,\pi\circ pr_1   \right)$ is an analytical continuation, with logarithmic singularities of $\left({\cal V}_2,g_2   \right)$; indeed
$
\displaystyle
 \pi\circ pr_1
\circ J=\pi\circ
\widetilde{\hbox{\j}}\circ g_2=g_2$;
but $\left({\cal V_2},g_2   \right)$ is connectible with $\left({\cal V},g     \right)$, hence
$\displaystyle \left(\Sigma,pr_2,J,\pi\circ pr_1   \right)$ is an analytical continuation of $\left({\cal V},g \right)$.

Thus there exists a continuous function $h:\Sigma\rightarrow S$,
holomorphic on the interior of $\Sigma$,
such that $\rho\circ h=pr_2$ : hence
$
\eta=pr_2\left(\xi,\eta   \right)=\rho\circ h
\left(\xi,\eta   \right)
\in\rho\left(S   \right)
$.
\QUAN
\section{Paths}

Let's start with a slight
reformulation of the notion of path: to achieve this goal, we adopt the point of view according to which a 'path' or even a 'curve' are analytical continuations of some initial germs, generally yielded by local solutions of systems of differential equations.

We shall also deal with the
{\it velocity field}
of a path:
to define it we shall need to single out
a vector field on its domain of definition,
which will have
to be related with the natural derivation field $d/dz$ on $\CI$.

Let $\M$ be a connected complex manifold: in the continuation, abusing language but following Wells (see e.g. \cite{wells} or \cite{gunros}), we shall name $T\M$ (resp.$T^{*}\M$) its holomorphic tangent (resp. cotangent) bundle and, more generally, ${\cal T}_{r}^{s}\M$ its holomorphic r-covariant and s-contravariant tensor bundle; as usual, $\Pi\colon{\cal T}_{r}^{s}\M\rightarrow\M$ will denote
 the natural projection.

A {\bf closed hypersurface} ${\cal F}$
in $\M$ is a closed subset such that there
exists
a maximal atlas $\{U_n\}$ for $\M$ and,
for each $n$, a holomorphic function
$\Psi_n$, not vanishing everywhere,
such that
$U_n\bigcap{\cal F}
=\{X\in U_n : \Psi_n(X)=0\}$.

The following definition is adapted from \cite{oneill}, definition 2.4 and lemma 2.5:
\begin{definition}
let ${\cal E}$ be a closed hypersurface in $\M$, $\N$ another connected complex manifold and $F\in{\cal O}(\M,\N)$:
{\bf an ${\cal E}$-meromorphic section} of
${\cal T}_{r}^{s}\N$ {\bf over} $F$ is a
holomorphic section $\Lambda$ of ${\cal T}_{r}^{s}\N$ over $F\vert_{\M\setminus
{\cal E}}$ such that
$\pi\circ\Lambda$
admits analytical
continuation up to
 the whole $\M$
and
for every $p\in{\cal E}$
and every coordinate system
$\left({\cal U},(z^1...z^n)
\right)$
around
$F(p)$, there exists a
neighbourhood $U$ of
$p$ and $r\cdot s$
pairs of $\CI-$valued
holomorphic functions
$\phi_{i_1...i_r},\
\psi_{l_1...l_s}$, with
 $\psi_{l_1...l_s}\not=0$
on
$U\setminus{\cal E}$, such that
$$
\Lambda
\left( dz^{l_1}...dz^{l_s},
\frac
{\partial}
{\partial z^{i_1}}...
\frac{\partial}
{\partial z^{i_r}}
\right)=
\frac{\phi_{i_1...i_r}}
{\psi_{l_1...l_s}}.
$$
\label{section}
\end{definition}

\begin{definition}
A {\bf path} in $\M$ is a
quadruple $Q_{\M}=
\left(S,\pi,j,F\right)$,
where
$S$ is a connected Riemann
 surface,
$\pi\colon S
\rightarrow\CI$
is a branched covering of
$S$ over $\pi(S)$,
$F:S\rightarrow\M$ is a
holomorphic mapping,
$U\subset{\CI}$ is an open
set wich admits a
holomorphic (hence open)
immersion $j\colon
U\rightarrow
S\setminus\Sigma$
such that $\pi\circ
j=id\vert_{U}$.
\label{path1}
\end{definition}

We are now turning to
define the {\bf velocity
field} of a path $Q_{\M}$:
it will be defined
as a suitable meromorphic
section over $F$ of the
holomorphic tangent bundle
$T\M$.
To achieve this purpose,
we need to lift the vector
field ${d}/{dz}$ on
$\CI$ with respect to $\pi$.

Of course, in general,
contravariant tensor fields
couldn't be lifted;
notwithstanding, we may get through this obstruction by keeping into account that $\CI$ and $S$ are one-dimensional manifolds and allowing the lifted vector field to be meromorphic: these matters are fathomed in next statements: recall that $P$ is the set of branch points of $\pi$.

\begin{lemma}
There exists a unique $P$-meromorphic vector field $\widetilde{{d}/{dz}}$ on $S$
such that, for every $r\in S\setminus P$, $\pi_{*}\vert_r\left(\widetilde{{d}/{dz}}\vert_{r}    \right)=\left({d}/{dz}\right)\vert_{\pi(r)}$.
\label{lift}
\end{lemma}
{\bf Proof:} Consider $\omega=\pi^{*}dz$ and  $\Lambda=\pi^{*}(dz\odot dz    )$ on $S$: the latter establishes an isomorphism between the holomorphic cotangent  and  tangent bundles
of $S\setminus P$.
Call $V$ the holomorphic vector field corresponding to $\omega$ in the above isomorphism: we claim that $V=\widetilde{{d}/{dz}}$ on $S\setminus P$.
To show this fact, we explicitely compute the components of $V$ with respect to
a maximal atlas ${\cal B}=\left\{(U_{\nu},\zeta_{\nu}    )\right\}$
for $S\setminus P$: let
$
\omega_{(\nu)\ 1}=
\omega(\partial/\partial \zeta_{(\nu)}    )$,
$
g_{(\nu)\ 11}=\Lambda(\partial/\partial \zeta_{(\nu)},\partial/\partial \zeta_{(\nu)}    )
$;
then, set ${V_{(\nu)}^1={\omega_{(\nu)\ 1}}/{g_{(\nu)\ 11}}}$ the collection $\left\{( {\cal U_{\nu}}, V_{(\nu)}^1  )\right\}$ of open sets and holomorphic functions is such that,
on overlapping local charts $(U_a,\zeta_a    )$ and $(U_b,\zeta_b   )$,
we have
$$
V_{(a)}^1=
\frac
{\omega_{(a)\ 1}}{g_{(a)\ 11}}
=\frac
{\omega_{(b)\ 1} (d\zeta_{(b)}/d\zeta_{(a)})}
{g_{(b)\ 11}({ d\zeta_{(b)}}/{d\zeta_{(a)}})^2}
=
{{V_{(b)}^1}
\frac
{d\zeta_{(a)}}{d\zeta_{(b)}}} ,
$$
that is to say that collection defines
a holomorphic
vector field. Now for every $r\in S\setminus P$,
$$
dz\vert_{\pi(r)}
(
\pi_{*}\vert_r
\widetilde{{d}/{dz}}\vert_{r}
)=\pi^*dz\vert_r
(
\widetilde{{d}/{dz}}\vert_{r}
)
%$
%$=
%\pi^*dz\vert_r
%(
%(
%\frac
%{1}
%{dz\vert_{\pi(r)}(\pi_*\partial/\partial\zeta\vert_r)}
%\frac{\partial}{\partial\zeta}
%\vert_r
%)
=\frac{\pi^*dz\vert_r (\partial/\partial\zeta\vert_r) }
{dz\vert_{\pi(r)}(\pi_*\partial/\partial\zeta\vert_r)
}
=1
,$$
hence
$
\pi_{*}\vert_r(\widetilde{{d}/{dz}}\vert_{r}    )=({d}/{dz})\vert_{\pi(r)}
%\iff dz\vert_{\pi(r)}(\pi_{*}\vert_r \widetilde{{d}/{dz}}\vert_{r}    )=1
$,
%We prove the right side to be true:
%$
%dz\vert_{\pi(r)}
%(
%pi_{*}\vert_r
%\widetilde{{d}/{dz}}\vert_{r}
%)=\pi^*dz\vert_r
%(
%\widetilde{{d}/{dz}}\vert_{r}
%)$
%$=\pi^*dz\vert_r
%(
%(
%\frac
%{1}
%{dz\vert_{\pi(r)}(\pi_*\partial/\partial\zeta\vert_r)}
%\frac{\partial}{\partial\zeta}
%\vert_r
%)
%=\frac{\pi^*dz\vert_r (\partial/\partial\zeta\vert_r) }
%{dz\vert_{\pi(r)}(\pi_*\partial/\partial\zeta\vert_r)
%}
%=1
%$,
%hence the left one is true too,
proving the asserted.

Let's prove that $\widetilde{d/dz}$ may be extended to a meromorphic vector field on $S$:
if $p\in P$ then
we can find local charts $(U,\psi    )$ around $p$, $(V,\phi)$ around $\pi(p)$,
and an integer $N>0$ such that $\phi\circ\pi\circ\psi^{-1}(u)=u^N$.
Now we have
$$
(\psi^{-1\ *}\pi^{*}\phi^{*}(dw)\frac d {du})(u)=
dw(\phi_{*}\pi_{*}\psi_{*}^{-1}\frac d {du})\vert_u)
=dw((\phi\pi\psi^{-1})^{\prime}\frac d {dw}))=Nu^{N-1};
$$
 but $\phi$ and $\psi$ are charts, hence $\pi^{*}dz $ itself is vanishing of order $N-1$ at $p$;
as already proved, $\pi_{*}\vert_r(\widetilde{{d}/{dz}}\vert_{r}    )=({d}/{dz})\vert_{\pi(r)}$ on $U\setminus \{p\}$ and, consequently,
$(\pi^{*}dz)(\widetilde{d/dz})=
dz(d/dz)=1$
on $U\setminus \{p\}$,
hence on $U$.
Now, in local coordinates,
$
(\pi^{*}dz)=\alpha d\phi$ and
$\widetilde{d/dz}=y\,{\partial}/{\partial\phi}
$,
where $\alpha$ is a suitable holomorphic function on $U$, vanishing of order $N-1$ at $p$ and $y$ is a holomorphic function on $U\setminus \{p\}$.
By the above argument, $y\alpha=1$, hence $y$ has a pole of order $N-1$ at $p$:
a similar argument holds for each isolated point in $P$,
proving the meromorphic behaviour of $\widetilde{d/dz}$.
\QUAN
\begin{definition}
A {\bf finite-velocity  point} of a path $Q_{\M}=\left(S,\pi,j,F,\M\right)$
is a point $r\in S$ such that $\widetilde {{d}/{dz}}$ is holomorphic at $r$.
\label{path3}
\end{definition}

We are ready to define the velocity field : let at first be $r$ a finite-velocity  point of $Q_{\M}$; since $\widetilde{{d}/{dz}}$ is holomorphic at $r$, we could define the holomorphic velocity at $r$ as $V_r=F_{*}\left((\widetilde{{d}/{dz}})\vert_r\right)$: now define
the mapping
$V\left(Q_{\M}\right)
\colon S\setminus
P \rightarrow
T\M$
 by setting
$
r\mapsto
\left(F,F_{*}
\left(\widetilde{{d}/{dz}}\
(r)\right)\right)\nonumber
$.

\begin{lemma}
The mapping
$V\left(Q_{\M}\right)$
can be extended to a P-meromorphic
section  of $T\M$ over $F$.
\label{velo_def}
\label{vf}
\end{lemma}
{\bf Proof:} Trivially $\Pi\circ V\vert_{R\setminus P}=F\vert_{R\setminus P}$.
Let's show the meromorphic behaviour of $V$: if $p\in P$ there is a neighbourhood $U$ of $p$ such that, for every local chart $\zeta\colon U\rightarrow\CI_w$ there exist holomorphic functions $f,g\in{\cal H}\left(\zeta(U)\right)$ such that
$$\widetilde{\frac{d}{dz}}\vert_{\zeta_{-1}(U)}=\zeta^{-1}_{*}\left(\frac{f}{g}(w)\frac{d}{dw}\vert_w\right);$$
moreover, for every local chart $\Psi=\left(u^1...u^m,du^1...du^m\right)$ in $T\M$ we obtain
\begin{eqnarray}
& &\Psi\circ V \circ \zeta^{-1}(w)=\nonumber\\
&=& \Psi\circ \left(F\circ \zeta^{-1}(w),F_{*}\vert_{\zeta^{-1}(w)}\left(\widetilde{\frac{d}{dz}}\vert_{\zeta^{-1}(w)} \right)\right) \nonumber\\
&=& \Psi\circ \left(
F\circ \zeta^{-1}(w),F_{*}\vert_{\zeta^{-1}(w)}\zeta^{-1}_{*}\left(\frac{f}{g}(w)\frac{d}{dw}\vert_w\right)\right) \nonumber\\
&=& \Psi\circ \left( F\circ \zeta^{-1}(w),\frac{f}{g}(w)\frac{d}{dw}(F\circ\zeta^{-1})(w)\right) \nonumber\\
&=& \Big(u^1\circ F\zeta^{-1}(w)...u^m\circ  F\zeta^{-1}(w),\nonumber\\
& &
\left.\frac{f}{g}(w)\frac{d}{dw}\left(u^1\circ F\circ\zeta^{-1}\right)(w)...\frac{f}{g}(w)\frac{d}{dw}\left(u^m\circ
F\circ\zeta^{-1}\right)(w))\right)  \nonumber
\end{eqnarray}
\QUAN

According to lemma \ref{velo_def},
the {\bf velocity field } of
a path $Q_{\M}=(S,\pi,j,F,\M)$
will be just the meromorphic mapping
$V\left(Q_{\M}\right)$.

\subsection{Definition of completeness}
\begin{definition}
\label{completezza}
A $\M$-valued path
 $\displaystyle
\left(S,\pi,j,F
\right)$ is {\TTT
complex-complete}
 provided that
$\CI\setminus
\displaystyle\pi\left(S
    \right)$ is a
finite set in the
complex plane;
a real-analytic
curve $\gamma$ in a
real-analytic manifold
$\N$ admitting a
complexification
$d:\N\rightarrow\M$ is
{\sl real-complete} (or, briefly,
{\sl complete})
provided that the Riemann surface
$(S,\pi,j,G)$ of $d\CIRC\gamma$ is such that
$\ERRE\setminus\pi(G^{-1}(d(\N)))$
is
a finite set.
\labelle{completessa}
\end{definition}

\section{Complex-Riemannian metric structures}
The intuitive geometry of the real euclidean space $\ERRE^3$ can be easily brought back to its natural inner product, which allows basic geometrical operations, like measuring the length of a tangent vector, or angles between tangent vectors: Riemannian real geometry generalizes all this to 'curved' spaces, which is based on the concept of positive definite bilinear forms: weakening definiteness to nondegeneracy leads us in the realm of Lorentz geometry, originating from the problems posed by Einstein's general
relativity theory.

A little bit less intuitive is the idea
of starting from the basic geometry of $\CI^3$
(meant as a 'complexification' of the usual real euclidean space) to get formal extension of the geometric properties of real 'curved' manifolds.
Introducing this complex environment could allow us to hope to get able to handle some types of metrical singularities which naturally arise in dealing with real manifolds with indefinite metrics.

It is immediately
seen that the nondegeneracy hypothesis
itself
should be dropped,
as the following considerations show
(see \cite{dubnovfom} p.186 ff):
consider the space
${\cal F}$ of antisymmetric covariant tensors of rank two
in Minkowski's space $\ERRE_{1,3}$:
electromagnetic fields are such ones.
Let $F\in{\cal F}$: we can write $F=
\sum_{i<j}F_{ij}dx^i\wedge dx^j$ where
$x^0...x^3$ are the natural coordinate functions
on $\ERRE_{1,3}$.
At each point, the space ${\cal F}_p$ of
all tensors in ${\cal F}$ evaluated at $p$ is a
six-dimensional real vector space; moreover, the
adjoint operator $*$ with respect to Minkowski's metric is such that $**=-1$:
all these facts imply that
${\cal F}_p$ could be thought of as a complex three dimensional
vector space ${\cal G}_p$ by setting $(a+\IM b)F=aF+b*F$.
Now $*$ is $SO(1,3)-$invariant, hence
$SO(1,3)$ is a group of (complex) linear transformations of ${\cal G}_p$,
preserving the quadratic form $\langle F,F \rangle=
-*\left(F\wedge (*F)+\IM  F\wedge F\right)$: this means
that this 'norm' is invariant by Lorentz transformations,
hence it is of relevant physical
interest.
If we introduce the following coordinate functions
on ${\cal G}_p$:
$
z^1=F_{01}-iF_{23}
$,
$
z^2=F_{02}+iF_{13}
$ and
$
z^1=F_{03}-iF_{12}
$,
we have that
%\begin{equation}
$
\displaystyle
\langle F,F \rangle=
(z^1)^2+(z^2)^2+(z^3)^2,
\label{complexeuclid}
$
%\end{equation}
 hence there naturally
arises the so called
{\it complex-Euclidean metric}
on $\CI^3$: on one hand, by changing coordinates
we are brought to a generic symmetric bilinear
form on $\CI^3$; on the other,
there arise 'poles'
if we attempt to extend
%(\ref{complexeuclid})
the above construction e.g. to $(\PI^1)^3$.
Now the idea of generalizing
to the curved framework is
quite natural:
let $\M$ be a complex manifold, ${\cal D}$ and ${\cal E}$ closed hypersurfaces in $\M$.

\begin{definition}
{\TTT A  holomorphic (resp.${\cal E}$-meromorphic) metric on} $\M$ is a holomorphic (resp.${\cal E}$-meromorphic)  section
$\Lambda\colon\M\rightarrow {\cal T}_{0}^{2}\M$ which is symmetric, that is to say, for every $m\in\M$ and every pair
of holomorphic tangent vectors
$V_m$ and $W_m$ at $m$,
there holds
$\Lambda(m)\left( V_m, W_m  \right)=\Lambda(m)\left( W_m, V_m  \right) ;
$.
The {\bf rank} of $\Lambda$
at $p\in\M$ is the rank of
the bilinear form $\Lambda(p)$;
$\Lambda$ is {\bf nondegenerate}
 at $p$ if $rk(\Lambda(p))=
dim(\M)$, {\bf degenerate}
otherwise;
if ${\cal D}$ is a
hypersurface in $\M$ and
$\Lambda$ is degenerate
only on $\cal D$, we shall
say that $\Lambda$ is
${\cal D}$-degenerate.
We say that $p$ is a
{\bf metrically ordinary point}
 in $\M$ if $\Lambda$ is
holomorphic and nondegenerate
at $p$.
\label{reg_point}
\label{rank_met}\label{hol_met}
\end{definition}
In the following we shall consider only metric which degenerate only on closed hypersurfaces.

\begin{definition}
A {\bf holomorphic
Riemannian manifold} is
a complex manifold endowed
with a holomorphic metric;
a {\bf nondegenerate holomorphic Riemannian manifold} is a complex manifold endowed with a nondegenerate holomorphic metric ;
a {\bf meromorphic  Riemannian manifold} is a complex manifold endowed with a
meromorphic metric.
\label{riemann}
\end{definition}
Thus, strictly speaking, all the above objects are pairs consisting in complex manifolds and  metrics, but we shall often understand metrics and denote them by the only underlying complex manifolds.

\subsection{The meromorphic Levi-Civita connexion}
We begin this section by introducing the holomorphic Levi-Civita connexion induced on a holomorphic nondegenerate Riemannian manifold by its metric structure: this is done in a quite similar way to that pursued in (real) differential geometry, apart from a slight difference, which naturally arises: the action of the Levi-Civita connexion is defined at first on 'local' vector fields, producing local ones as well, then it is globalized as a collection of local operators.

Let now $\left(\M,\Lambda    \right)$ be a nondegenerate Riemannian holomorphic manifold, ${\cal A}$ a maximal atlas for $\M$, ${\cal U}\in{\cal A}$ a domain of a local chart.
Let also
${\cal X}\left({\cal U}
 \right)$ be the Lie
algebra of holomorphic
vector fields on ${\cal U}$
and
${\cal O}\left({\cal U}
    \right)$ the ring
of holomorphic functions
 on ${\cal U}$.

\begin{definition}
A {\bf connexion on ${\cal U}$} is a mapping
$D \colon{\cal X}\left({\cal U}    \right)\times{\cal X}\left({\cal U}    \right)\rightarrow {\cal X}\left({\cal U}    \right)$ such that:
{\bf (D1)}
$D_V W$ is ${\cal H}\left({\cal U}    \right)$-linear in $V$;
{\bf (D2)}
$D_V W$ is $\CI$-linear in $W$ and
{\bf (D3)}
$D_V\left(fW    \right)=\left(Vf\right)W+fD_VW$ for every $f\in{\cal H}\left({\cal U}    \right)$.
\label{met_conn_germ}
\end{definition}

$D_V W$ is called the {\bf covariant derivative } of $W$ with respect to $V$
 in the connexion $D$. By axiom (D1), $D_V W$ has tensor character in $V$, while axiom (D3) tells us that it is not a tensor in $W$.

Our next step is to show that there is a unique connexion characterized by two further properties, (D4) and (D5) below, namely being anti-Leibnitz like with respect to the Lie bracket operation and Leibnitz like with respect to the metric.
In the following we use the alternative notation $\left\langle V,W   \right\rangle$ instead of $\Lambda\left(V,W    \right)$.
\begin{lemma}
Let ${\cal U}$ be an open set belonging to a maximal atlas ${\cal A}$ for the nondegenerate holomorphic Riemannian manifold $\M$. If $V\in{\cal X}\left({\cal U}    \right)$, let $V^*$ be the holomorphic one-form on ${\cal U}$ such that $V^*(X)=\left\langle V,X   \right\rangle$ for every $X\in{\cal X}\left({\cal U}    \right)$: then the mapping $V\mapsto V^*$ is a ${\cal O}$-linear isomorphism from ${\cal X}\left({\cal U}    \right)$ to ${\cal X}^*\left({\cal U}    \right)$.
\label{vefo}
\end{lemma}
{\bf Proof:} Since $V^*$ is ${\cal O}$-linear, it is in fact a one-form, and $V\mapsto V^*$ is ${\cal O}$-linear too.
We claim:
\begin{itemize}
\item[(a)] if $\left\langle V,X   \right\rangle=\left\langle W,X   \right\rangle$ for every $X\in{\cal X}\left({\cal U}    \right)$ then $V=W$;
\item[(b)] given any one-form $\omega\in{\cal X}^*\left({\cal U}    \right)$ there is a uique vector field $V\in{\cal X}\left({\cal U}    \right)$ such that $\omega(X)=\left\langle V,X   \right\rangle$ for every $X\in{\cal X}\left({\cal U}    \right)$.
\end{itemize}

Let $U=V-W$;
the nondegeneracy
of the metric tensor
implies that,
if $p\in{\cal U}$
and
$
\left\langle U_p,X_p
   \right\rangle=0$
for every $X\in{\cal X}
\left({\cal U}    \right)$,
then
 $U=0$;
this proves (a).

To prove (b), let $\left(z^1...z^N    \right)$ be local coordinates on ${\cal U}$.

Then $\omega=\sum_{i=1}^N\omega_i dz^i$; let $\{g_{ij}\}$ be the representative matrix of $\Lambda\vert_{\cal U}$ in $\left(z^1...z^N    \right)$: by nondegeneracy,
it admits a holomorphic
inverse matrix $\{g^{ij}\}$:
set now $V=
\sum_{j=1}^N\left
(\sum_{i=1}^N g^{ij}\omega_i
   \right)
\frac{\partial}{\partial z^j}
$.

We have
$
\left\langle V,X
  \right\rangle =
 \left\langle
 \sum_{j=1}^N
\left(\sum_{i=1}^N g^{ij}
\omega_i
\right)\frac{\partial}
{\partial z^j},
\sum_{k=1}^N X^k
\frac{\partial}{\partial z^k}
  \right\rangle=
\sum_{ijk}g^{ij}\omega_i
X^k g_{jk}
= \sum_{ik}\delta_k^i X^k
\omega^i
=\sum_k X^k\omega^k
=\omega\left(X    \right)
$.
\QUAN

The following theorem can be proved
exactly as in classical differential
geometry.
\begin{theorem}
Let ${\cal U}$ be an open set belonging to a maximal atlas ${\cal A}$ for the nondegenerate holomorphic Riemannian manifold $\M$. There exists a unique connexion $D$ on ${\cal U}$, called the Levi-Civita connexion, such that:
\begin{itemize}
\item[(D4)]
$\left[V,W\right]=D_V W-D_W V$;
\item[(D5)]
$X\left\langle V,W   \right\rangle=\left\langle D_X V,W   \right\rangle + \left\langle V,D_X  W   \right\rangle $ for every $X,V,W\in {\cal X}\left({\cal U}    \right)$.
\end{itemize}
Moreover $D$ is characterized by the 'Koszul's formula':
$
2\left\langle D_V W,X   \right\rangle = V\left\langle W,X   \right\rangle+
W\left\langle X,V   \right\rangle - X\left\langle V,W   \right\rangle
-\left\langle V, [W,X]   \right\rangle
+\left\langle W, [X,V]   \right\rangle
+\left\langle X, [V,W]   \right\rangle
$,
for every $X,V,W\in {\cal X}\left({\cal U}    \right)$.
\label{levicivita}

\end{theorem}

If we have to emphasize the open set ${\cal U}$ in theorem \ref{levicivita} we shall write $D\left[ {\cal U}   \right]$ instead of $D$:  if $\ {\cal U}_1,{\cal U}_2 \subset \M$ in a maximal atlas
${\cal A}$ for $\M$ are
overlapping open sets,
then ${\cal U}_{1}
\bigcap{\cal U}_{2} $
is in ${\cal A}$
too and
$
D\left[{\cal U}_{1}    \right]\vert_{{\cal X}\left({\cal U}_{1}\bigcap{\cal U}_{2}    \right)}=
D\left[{\cal U}_{1}    \right]\vert_{{\cal X}\left({\cal U}_{1}\bigcap{\cal U}_{2}    \right)}
$,
hence we can collect all
local definitions of
Levi-Civita connexions:
\begin{definition}
the {\bf Levi-Civita connexion}
 (or {\bf metric connexion})
$D$ of
$\left(\M,\Lambda    \right)$
is the collection consisting
of all the metric connexions
$\{D\left[{\cal U}_i
    \right]\}_
{i\in I}$ as
$
{\cal U}_i $ runs
over any maximal atlas
${\cal A}=\left(
\{
{\cal U}_i
\}
\right)_{i\in I}$ on $\M$.
\label{met_conn}
\end{definition}

So far we have studied
nondegenerate
holomorphic Riemannian
manifolds: this situation
is quite similar to real
Riemannian geometry.

Things are different,
instead, if we allow metrics
to have meromorphic behaviour, or to lower somwhere in their ranks.
These metric 'singularities'
 will be generally supposed
to lie in closed
hypersurfaces;
Levi Civita  connexions
may still be defined, but,
as one could expect, they
will turn out to be
themselves 'meromorphic'.

Let now
$\left(\N,\Lambda    \right)$
be a meromorphic
Riemannian manifold
admitting closed
hypersurfaces ${\cal D}$
and ${\cal E}$ such that
$\Lambda
\vert_{\N\setminus{\cal E}}$
 is holomorphic and
$\Lambda
\vert_{(\N\setminus{\cal E})
\setminus{\cal D}}$ is
nondegenerate.
Since $\N\setminus{\cal E}$
is connected, we have that
$(\N\setminus{\cal E})
\setminus{\cal D},\Lambda
\vert_{(\N\setminus{\cal E})
\setminus{\cal D}}$ is a
nondegenerate holomorphic
Riemannian manifold
admitting, as such,
a canonical holomorphic
Levi-Civita connexion $D$.

Now,
if $p\in{\cal D}
\bigcup{\cal E}$
and $V,W$ are
holomorphic vector
fields in a neighbourhood
${\cal V} $ of $p$,
it will result that
we are able to define
the vector field $D_V W$ on
${\cal V}
\setminus
\left({\cal D}\bigcup{\cal E}
    \right)$,
and this will be
a meromorphic vector field.

Let's state all this more
precisely:
\begin{definition}
Let $Z=(z^1\cdots z^m)$
be
a coordinate system
on
an open set ${\cal U}
\subset\N$:
the {\bf Christoffel symbols }
of $Z$
 are those complex valued
functions, defined on
${\cal U}\setminus\left
({\cal D}\bigcup{\cal E}
    \right)$ by setting
$\Gamma_{ij}^k =
dz^k\left
(  D_{\frac{\partial}
{\partial z^i}}\left
({\frac{\partial}
{\partial z^j}}   \right)
 \right)$.
\label{christoffel}
\end{definition}
Now the representative matrix
$(g_{ij})$
of $\Lambda$ with respect
to the coordinate system $Z$
is holomorphic in ${\cal U}$,
 with nonvanishing determinant
 function on
${\cal U}\setminus
\left({\cal D}\bigcup{\cal E}    \right)$; as such it admits a inverse
matrix ${g^{ij}}$,
whose coefficients hence
result in being
${\cal D}\bigcup{\cal E}
$-meromorphic functions.

\begin{lemma}
{\bf (a)}
$D_{\frac{\partial}
{\partial z^i}}
\left(\sum_{j=1}^m W^j
{\frac{\partial}{\partial z^j}
 }   \right)=
\sum_{k=1}^m
\left(\frac{\partial W^k}
{\partial z^i}+\sum_{j=1}^m
\Gamma_{ij}^k W^j
  \right)\frac{\partial}
{\partial z^k}$
as meromorphic vector fields;
{\bf (b)}
$2\Gamma_{ij}^k
=\sum_{m=1}^N g^{km}
\left(-g_{ij,m}+g_{im,j}
+g_{jm,i}\right)=2\Gamma
_{ij}^k $
as meromorphic functions.
\label{christoffel_2}
\end{lemma}

{\bf Proof:} At first note that the operation of associating Christoffel symbols
to a coordinate system is compatible with restrictions, in the sense that the Christoffel symbols of the restriction of $Z$ to a smaller open set are its Christoffel symbols restricted to  that  set.
Now, if $
p\in
{\cal U}\bigcap
\{n\in\N:\Lambda
\hbox
{ is holomorphic
and nondegenerate at $n$}\}
$
and ${\cal V}_p\subset {\cal U}$ is a neighbourhood of $p$, contained in ${\cal U}$, we have that $\Lambda$ is holomorphic and nondegenerate in ${\cal V}_p$: hence
(a):  by Koszul's formula
 we have
$$
2\sum_{a=1}^N \Gamma_{ij}^a g_{am}=2\left\langle D_{\frac{\partial}{\partial z^i}} \frac{\partial}{\partial z^j},\frac{\partial}{\partial z^m}  \right\rangle=
\frac{\partial}{\partial z^i}g_{jm}+\frac{\partial}{\partial z^j}g_{im}+\frac{\partial}{\partial z^m}g_{ij};
$$
multiplying both side by $g^{mk}$ and summing over $m$ yields the desired result; (b) follows immediately from (D3) of definition \ref{met_conn_germ}.
Now the fact that (a) and (b) hold in fact on ${\cal U}$ follows by analytical continuation: note that this result does not depend on the choice of $p$.
\QUAN
\begin{proposition}
For every pair $V,W$
of holomorphic vector
fields on the open set
${\cal U}$ (
belonging to a maximal atlas)
in the meromorphic Riemannian
manifold $\left(\N,\Lambda
     \right)$, $D_V W$ is a
well defined vector field,
holomorphic on
${\cal U}\bigcap
\{n\in\N:\Lambda\hbox
{ is holomorphic and
nondegenerate at $n$} \}$ and may be extended to a meromorphic vector field on ${\cal U}$.
\end{proposition}
{\bf Proof:} There exist holomorphic functions $\{V^i\}$, $\{W^j\}$ and a coordinate system $Z=\left(z^1.....z^N    \right)$
on ${\cal U}$ such that
$
\displaystyle{V=\sum_{i=1}^N V^i\frac{\partial}{\partial z^i}}
$
and
$
\displaystyle{W=\sum_{j=1}^N W^i\frac{\partial}{\partial z^j}}
$.
By lemma \ref{christoffel_2}(a),
$$
D_V W=\sum_{i=1}^N V^i D_{\frac{\partial}{\partial z^i}}\left(\sum_{j=1}^N W^j   \frac{\partial}{\partial z^i} \right)
=\sum_{k=1}^N \left(\sum_{i,j=1}^N V^i\left( \frac{\partial W^k}{\partial z^i}+\Gamma^k_{ij}W^j   \right)    \right)\frac{\partial}{\partial z^k}:
$$
this is a vector field whose components are meromorphic functions.
\QUAN
Summing up, we yield:
\begin{definition}
Given a ${\cal D}$-degenerate and ${\cal E}$-meromorphic Riemannian manifold $\left(\N,\Lambda    \right)$, with ${\cal
D} $ and ${\cal E}$ closed hypersurfaces in $\N$, the {\bf Levi-Civita metric connexion} (or {\bf meromorphic metric
connexion}) of $\N$ is the collection consisting of the metric connexions $\{D\left[{\cal U}_i\setminus({\cal D}\bigcup
{\cal E}  ) \right]\}_{i\in I}$ as ${\cal U}\}_i $ runs over any maximal atlas ${\cal B}=\left(\{{\cal U}\}_i    \right)_{i\in I}$ on $\N$.

\end{definition}

\subsection{Meromorphic parallel translation}
We turn now to study vector fields on paths: an obvious example is the velocity field
(
see lemma
\ref{vf}): just as in semi-Riemannian geometry, there is a natural way of defining the rate of change $X^{\prime}$ of a meromorphic vector field $X$ on a path.
We study at first paths
 with values in a nondegenerate holomorphic Riemannian manifold $\M$:
let
$Q_{\M}=\left
(S,\pi,j,\gamma,\M
\right)$ be a path in $\M$,
$P$ be the set of branch
points of $\pi$,
$r\in S\setminus P$ be
 a finite-velocity point
of $Q_{\M}$.
Moreover, let
${\cal V}\subset
S\setminus P$ be a
neighbourhood of $r$
such that $\gamma\left
({\cal V}\right)$ is
contained in a local
 chart in $\M$,
${\cal H\left(V
    \right)}$ be
the ring of holomorphic
functions on ${\cal V}$,
${\cal X}_{\gamma}\left
({\cal V}    \right)$ the
 Lie algebra of holomorphic
 vector fields over $\gamma$
on ${\cal V}$.

Due to the locally nondegenerate holomorphic
environment, the following proposition
can be proved in  quite a classical fashion.
\begin{proposition}
There exists a unique mapping
$\nabla_{\gamma^{\prime}}
\colon{\cal X}_{\gamma}
\left({\cal V}    \right)
\rightarrow
{\cal X}_{\gamma}
\left({\cal V}    \right)
$,
called {\bf induced covariant derivative} on $Q_{\M}$ in ${\cal V}$, (or on $\gamma$ in ${\cal V}$) such that:
\begin{eqnarray}
\hbox{(a)}\ & & \nabla_{\gamma^{\prime}}\left(aZ_1+bZ_2    \right)=a\nabla_{\gamma^{\prime}}
Z_1+b\nabla_{\gamma^{\prime}} Z_2,\quad a,b\in\CI\nonumber;\\
\hbox{(b)}\  & & \nabla_{\gamma^{\prime}} \left(hZ\right)=\left(\widetilde{\frac{d}{dz}}h    \right)Z
+h\nabla_{\gamma^{\prime}} Z,\quad h\in{\cal H\left(V    \right)}; \nonumber\\
\hbox{(c)}\ & & \nabla_{\gamma^{\prime}} \left(V\circ\gamma \right)(r)=
D_{\gamma_{*}\vert_r(\widetilde{\frac{d}{dz}}\vert_r)}\ r\in{\cal V},\nonumber
\end{eqnarray}
where $V$
is a holomorphic vector field in a neighbourhood of $\gamma(r)$.
Moreover,
$$\widetilde{\frac{d}{dz}}\left\langle X,Y\right\rangle=\left\langle \nabla_{\gamma^{\prime}}
 X,Y   \right\rangle+\left\langle X,\nabla_{\gamma^{\prime}} Y \right\rangle
\quad X,Y\in {\cal X}_{\gamma}({\cal V}).$$
\label{traspar}
\end{proposition}

Now let ${\cal R}=\{{\cal V}_k\}_{k\in K}$ be a maximal atlas for $S\setminus P$; we may assume that, for every $k$,
maybe shrinking ${\cal V}_k   $, $\gamma\left( {\cal V}_k   \right)$ is contained in some local chart ${\cal U}_i$ in the
already introduced atlas ${\cal A}$ for $\M$.

By proposition \ref{traspar}, if ${\cal V}_1  $ and ${\cal V}_2  $ are overlapping open sets in ${\cal R}$, ${\cal V}_1 \bigcap{\cal V}_2\in {\cal R} $ too, and
$
\nabla_{\gamma^{\prime}}
\left[ {\cal V}_1   \right]
\vert_{{\cal V}_1 \bigcap{\cal V}_2}
=\nabla_{\gamma^{\prime}}
\left[ {\cal V}_2   \right]
\vert_{{\cal V}_1 \bigcap{\cal V}_2}
$.

Now let's complete ${\cal R}$ to an atlas ${\cal S}$ for $S$: keeping into account that the local coordinate expression of the induced covariant derivative is
$$
\nabla_{\gamma^{\prime}}Z=\sum_{k=1}^m\left( \widetilde{\frac{d}{dz}}Z^k+\sum_{i,j=1}^m\Gamma_{ij}^k
\widetilde{\frac{d}{dz}}\left(u^i\circ\gamma\right)Z^j   \right)\frac{\partial}{\partial u^k}.
$$
and arguing in the same way as about the meromorphic Levi-Civita connexion, we are able to show that pairs of holomorphic vector fields on $\gamma$ are transormed into $P$-meromorphic vector fields on $\gamma$.

\begin{definition}
The $P$-meromorphic {\bf induced covariant derivative}, or the $P$-meromorphic parallel translation on a path $Q_{\M}=\left(S,\pi,j,\gamma,\right)$ with set of branch points $P$ and taking values in a nondegenerate Riemannian manifold $\M$ is the collection consisting of the induced covariant derivatives $\nabla_{\gamma^{\prime}}\left[ {\cal V}_k \setminus P  \right] $ as ${\cal V}_k $ runs over a maximal atlas ${\cal S}=\left(\{{\cal V}_k\}    \right)_{k\in K}$ on $S$.
\label{indcovdev}
\end{definition}

Let's turn now to dealing
with meromorphic parallel
translations induced on a
path $Q_{\N}=\left(T,\varrho,
j,\delta\right)$,
in a meromorphic Riemannian manifold
$(\N,\Lambda)$ admitting closed hypersurfaces ${\cal D}$ and ${\cal E}$ such that
$\Lambda\vert_{\N\setminus
{\cal E}}$ is holomorphic and
$\Lambda
\vert_{(\N\setminus{\cal E})
\setminus{\cal D}}$
is nondegenerate.
We set ${\cal F}={\cal D}
\bigcup{\cal E}$ and
restrict our attention to
paths starting at
metrically ordinary points.
\begin{lemma}
Set $\M=\N\setminus{\cal F}$, $S=\delta^{-1}(\M)$: then $T\setminus S$ is discrete, hence $S$ is a connected Riemann surface.
\end{lemma}
{\bf Proof:} Suppose that there exists a subset ${\cal V}\subset T\setminus S$ admitting an accumulation point $t\in{\cal V}$ and consider a countable atlas for ${\cal B}=\{U_n\}_{n\in\ENNE}$ for $\N$ such that, for every $n$, there exists $\Psi_n\in{\cal O}\left( \{U_n\}   \right)$ such that
$
U_n\bigcap{\cal F}=\{X\in U_n : \Psi_n=0\}
$.

Set $\delta^{-1}(U_n)=T_n\subset T$ and suppose, without loss of generality, that $\delta(t)\in U_0$.

We have $\Psi_0\circ\delta\vert_{{\cal V}\cap T_0}=0$ and $t\in{\cal V}\cap T_0$ is an accumulation point of ${\cal V}\cap T_0$ , hence $\Psi_0\circ\delta\vert_{ T_0}=0$ and $\delta(T_0)\subset {\cal F}$.

Suppose now that $T_N\not=\emptyset$ for some $N$: we claim that this implies $\delta(T_N)\subset{\cal F}$: to prove the asserted, pick two points $\tau_0\in T_0$ and $\tau_n\in T_n$ and two neighbourhoods $T^{\prime}_0$, $T^{\prime}_N$ of $\tau_0 $ and $\tau_n$ in $T_0$ and $T_n$ respectively, such that $\varrho\vert_{T^{\prime}_0}$ and $\varrho\vert_{T^{\prime}_N}$ are biholomorphic functions.
Now the function elements $\left(\varrho(T^{\prime}_0),\delta\circ\left( \varrho\vert_{T^{\prime}_0}   \right)^{-1}    \right)$ and $\left(\varrho(T^{\prime}_N),\delta\circ\left( \varrho\vert_{T^{\prime}_N}   \right)^{-1}    \right)$ are connectible, hence there exists a finite chain $\{W_{\nu}\}_{\nu=0...L}$ such that $W_0=\varrho(T^{\prime}_0)$, $W_L=\varrho(T^{\prime}_N)$, $W_{\nu}\bigcap W_{\nu+1}\not=0$ for every $\nu$.

Without loss of generality,
we may suppose that each
$W_{\nu}$ admits a
holomorphic, hence open,
immersion
$j_{\nu}\rightarrow T$,
hence, setting
$
S_0=T_0$,
$
S_{\lambda}=
j_{\lambda}(W_{\lambda})$
 for
$
\lambda=1...L
$ and
$
S_{L+1}=T_N
$
yields a finite chain of open subsets $\{S_{\lambda}\}_{\lambda=0...M}$ of $T$ connecting $T_0$ and $T_N$.

Let's prove, by induction, that, for every $\lambda$, $\delta(S_{\lambda})\subset{\cal F}$.

$\bullet$ {\bf At first}
recall that
$\delta(S_0)\subset
U_0\bigcap{\cal F}$
as already proved;
suppose now that
$\delta(S_{k-1})
\subset {\cal F}$.
We have $S_{k-1}
\bigcap S_k\not=\emptyset$,
hence $\delta(S_{k-1})
\bigcap\delta(S_{k})
\not=\emptyset$.

For every $m$ set
$\Sigma_{km}=
\delta(S_{k-1})
\bigcap\delta(S_{k})
\bigcap U_m
$:
 if $\Sigma_{km}\not=\emptyset$,
 then
$\Psi_m\circ\delta\equiv 0$ on
${\delta^{-1}
(\Sigma_{km})\bigcap
S_{k-1}\bigcap S_k }$;
 but $\delta^{-1}
(\Sigma_{km})\bigcap
S_{k-1}\bigcap S_k$
is open in
$\delta^{-1}\left(
\delta(S_k)\bigcap U_m
   \right)\bigcap S_k$,
thus
$\Psi_m\circ\delta\equiv 0$
on
$
{\delta^{-1}
\left(\delta(S_k)\bigcap
U_m    \right)\bigcap S_k}
$, that is to say $\delta(S_k)\bigcap U_m\subset{\cal F}$.

$\bullet$ {\bf On the other hand},
if $\Sigma_{km}=\emptyset$, but $\delta(S_k)\bigcup U_m\not=\emptyset$ we claim that $\delta(S_k)\bigcap U_m\subset{\cal F}$ as well: proving this requires a further induction: pick a $U_M$ such that $\Sigma_{kM}\not=\emptyset$ and a finite chain of open sets
${\cal B}^{\prime}=\{U^{\prime}_{\mu}\}_{\mu=0...J}\subset {\cal B}$ (with $U_{\mu}^{\prime}\bigcap\delta(S_k)\not=\emptyset$ for each $\mu$) connecting $U_M$ and $U_m$.
Since $\Sigma_{kM}\not=
\emptyset$,
$\delta(S_k)\bigcap
 U_0^{\prime}=\delta(S_k)
\bigcap U_M\subset {\cal F}$.

Suppose by induction that
$\delta(S_k)\bigcap
 U_{l-1}^{\prime}\subset
{\cal F}$.

Then
$\Psi_l\circ\delta
\equiv 0$
on ${\delta^{-1}
\left(\delta(S_k)\cap
U_{l-1}^{\prime}\cap
U_{l}^{\prime}   \right)
\cap S_k}$,
hence
$
\Psi_l\circ\delta\equiv 0$
on
$
{\delta^{-1}\left(\delta(S_k)\cap U_{l}^{\prime}   \right)\cap S_k}
$
i.e. $\delta(S_k)\bigcap U_{l}^{\prime}  \subset {\cal F}$: this ends the induction and eventually implies
$\delta(S_k)\bigcap U_{m}= \delta(S_k)\bigcap U_{J}^{\prime}\subset {\cal F}$.

Summing up,
$\delta(S_k)=\bigcup_m
 \left( \delta(S_k)\bigcap U_m
  \right)\subset {\cal F}$,
 for each $k$.
Hence
$\delta(T_N)=\delta(S_M)
\subset{\cal F}$
and eventually
 $\delta(T)=\delta\left(\bigcup_{N\in\ENNE} T_N    \right)\subset {\cal F}$,
 hence $\delta$ cannot start at a point in $\N\setminus{\cal F}$.
\QUAN

In the following considerations, there will still hold all notations introduced in
preceding lemma:
given a path
$Q_{\N}=\left(T,\varrho,j,
\delta\right)$,
set
 $\pi=\varrho\vert_S$,
$\gamma=\delta\vert_S$ and note that, since $Q_{\N}$ is starting from a metrically ordinary point $m$, $j$ may be
supposed to take values in fact in $S$; since the preceding lemma shows that $S$ is a connected Riemann surface,
$Q_{\M}=\left(S,\pi,j,\delta\vert_S\right)$ is in fact a path in $\M$, which we call the {\bf depolarization} of $Q_{\N}$. But $\M$ is a nondegenerate holomorphic Riemannian manifold, hence if $P$ is the set of branch points of $\pi$, there is a $P$-meromorphic induced parallel translation on $Q_{\M}$, got following definition \ref{indcovdev} and its substratum.
Finally, we introduce a maximal atlas ${\cal T}$ for $T$ and yield the following:
\begin{definition}
Let $\left(\N,\Lambda    \right)$ be a ${\cal E}$- meromorphic and ${\cal D}$-degenerate Riemannian manifold,
$\M=\N\setminus\left({\cal D}\bigcup{\cal E}    \right)$,
$Q_{\N}=\left(T\bigcup,\varrho
,j,\delta
\right)$ a path:
the {\bf $\left(P\bigcup\delta^{-1}\left( {\cal D}\bigcup{\cal E}   \right)    \right)$-meromorphic induced covariant derivative} on $Q_{\N}$
 is the collection consisting of all induced covariant derivatives $\nabla_{\gamma^{\prime}}\left[ {\cal V}_k \bigcap S
 \right] $ as ${\cal V}_k$ runs over a maximal atlas ${\cal T}=\left(\{{\cal V}_k\}    \right)_{k\in K}$ for  $T$ and
 $Q_{\M}=\left(S,\pi,j,\delta\vert_S    \right)$ is the depolarization of $Q_{\N}$.
\label{indcovdev2}
\end{definition}
\subsection{Geodesics}
\begin{definition}
A meromorphic (in particular, holomorphic) vector field $Z$ on a path
$Q_{\M}=\left(S,\pi,j,\gamma\right)$
is {\bf parallel} provided that $\nabla Z=0$ (as a {\bf meromorphic} field on $Q_{\M}$).
\label{parall}
\end{definition}
\begin{definition}
The {\bf acceleration} $\aleph\left( Q_{\M}\right)$ of $Q_{\M}$ is the meromorphic field $\nabla \left(V\left(Q_{\M}\right)\right)$ on $Q_{\M}$ yielded by the induced covariant
derivative of its velocity field;
\label{acceleration}
the {\bf speed} of a path is the 'amplitude' function of its velocity field:
 $S\left(Q_{\M}    \right)(r)=
 \left\langle \gamma_{*}
\vert_r\left(\widetilde{\frac{d}{dz}}\right), \gamma_{*}\vert_r\left(\widetilde{\frac{d}{dz}}   \right)   \right\rangle$. This is a meromorphic function.
A path is {\bf null} provided that its speed is zero everywhere. \label{speed}
\end{definition}
\begin{definition}
A {\bf geodesic} in a meromorphic (in particular, holomorphic) Riemannian manifold is a path whose
velocity field is parallel, or, equivalently, one of zero acceleration (see definition \ref{acceleration}).
A geodesic is {\bf null} provided that so is as a path.
\label{geodesic}
\end{definition}
The local equations of elements of geodesics $\left(U,\beta   \right)$
$
\BETA^{\bullet\bullet}{}^{k}+\sum_{i,j=1}^N \Gamma_{ij}^k(\beta)\BETA^{\bullet}{}^{i}\BETA^{\bullet}{}^{j}=0$,
$
(k=1.....N)
$
are a system of $N$  second-order ordinary differential equations in the complex domain, with meromorphic coefficients, in turn equivalent to an autonomous system of $2N$ first-order equations, hence, as a consequence of the general theory (see theorem \ref{gerexuq}) we have the following
\begin{theorem}
For every metrically ordinary point $p\in\M$, every holomorphic tangent vector $V_p\in T_p\M$ and every $z_0\in\CI$, there exists a unique germ
\hbox{\boldmath{}$\beta_{z_0}$\unboldmath}
of geodesic such that
$\hbox{\boldmath{}$\beta_{z_0}$\unboldmath $(z_0)=p$} $
and
$
\hbox{\boldmath{}$\beta_{z_0\ *}$\unboldmath$(d/dz)\vert_{z_0}=V_p$}
$;
moreover any analytical continuation of \hbox{\boldmath{}$\beta_{z_0}$\unboldmath}
is a geodesic.
\end{theorem}

\section{Completeness theorems}
\subsection{Complex
warped products}
In this section we shall be concerned with warped products of Riemann surfaces, each one endowed with some meromorphic metric: in this framework we shall prove a geodesic completeness criterion.

Let now $\,{\cal U}_i$,
$\left(i=1....N   \right)$, $N\geq 2$ be either a copy of the unit ball in the complex plane, or the complex plane itself, whose coordinate function we shall call $u^i$.

Moreover, let each $\,{\cal U}_i\,$ be endowed with a (not everywhere vanishing) meromorphic metric, which we denote by $b_1(u^1)\,du^1\odot du^1$ on $\,{\cal U}_1\,$, or by $f_i(u^i)\,du^i\odot du^i$ if $i\geq 2$, where
both $b_1$ and the $f_i$'s are nonzero meromorphic functions.

Consider now the meromorphic Riemannian manifold

$$
{\cal U}={\cal U}_1\times_{a_2(u^1)}
{\cal U}_2\times_{a_3(u^1)}
{\cal U}_3\times
........
\times_{a_N(u^1)}{\cal U}_N,
$$
where the $a_k$'s ($k\geq 2$) are
nonzero meromorphic warping functions defined on ${\cal U}_1$, i.e. depending solely on $u^1$.

We could write down the meromorphic metric $\Lambda$ of ${\cal U}$ in the form

$$
\Lambda\left(u^1.....u^N\right)=
b_1(u^1)\,du^i\odot du^i+
\sum_{i=2}^N a_i(u^i)f_i(u^i)
\,du^i\odot du^i.
$$

In other words, the matrix of $\Lambda $, with respect to the canonical coordinates of ${\cal U}$, inherited
from $\CI^N$, is of the form
$
(g_{ik})=diag\left(
b_1(u^1),
a_2(u^1)f_2(u^2),
a_3(u^1)f_3(u^3),
...
a_N(u^1)f_N(u^N)
\right)
$.
The following lemma can be proved by easy calculations:
\begin{lemma}
\label{connessione}
The meromorphic Levi-Civita
connexion
induced on ${\cal U}$
by ${\Lambda}$
admits the following
Christoffel symbols:
$
\displaystyle
2\Gamma_{11}^1=
{b_1^{\prime}(u^1)}/{b_1(u^1)}$;
$
\displaystyle
\Gamma_{ij}^1=0\ \hbox{if\
  } i\not=j$;
$
\displaystyle
2\Gamma_{ii}^1=
-\left[{a_i^{\prime}
(u^1)f_i(u^i)}    \right]/
{b_1(u^1)}\ \ \hbox{if\  }
1\leq i\leq N
$;
$
\displaystyle
2\Gamma_{kk}^k
={f_k^{\prime}(u^k)}/
{f(u^k)}$
if $2\leq k\leq N$
and
$
\displaystyle
2\Gamma_{ik}^k =
{a_k^{\prime}(u^1)}/
{a_k(u^1)}$
if $i=1$ and $2\leq k\leq N$.
Finally,
$
\displaystyle\
\Gamma_{ij}^k=0$
otherwise.
\end{lemma}
As an immediate consequence, we have:
\begin{lemma}
\label{equazionisecondoordine}
Each element of geodesic of $\left({\cal U},\Lambda   \right)$ satisfies the following system of $N$ ordinary differential equations in the complex domain:
\end{lemma}
\begin{equation}
\label{equazioninormali}
\cases{
\U^{\bullet\bullet}
{}^1+
\frac
{b_1^{\prime}(u^1)}
{2b_1(u^1)}
\left(\U^{\bullet}
{}^1\right)^2
-\sum_{l=2}^N
\frac{a_l^{\prime}
(u^1)f_l(u^l)}
{2b_1(u^1)}
\left(\U^{\bullet}{}^l
\right)^2=0\cr
\U^{\bullet\bullet}{}^k+
\frac{f_k^{\prime}(u^k)}
{2f_k(u^k)}
\left(\U^{\bullet}{}^k
\right)^2
+\frac{a_k^{\prime}
(u^1)}{a_k(u^1)}
\left(\U^{\bullet}{}^1
\right)\left(\U^{\bullet}
{}^k \right)
=0,\ k=2...N,
}
\end{equation}
provided that
it starts at a metrically
ordinary point.
Here, and in the following,
$u^k=u^k(z)$.

\begin{lemma}
\label{integraleprimo1}
The system
(\ref{equazioninormali})
 of the differential equations of elements of geodesics
$\displaystyle
z\mapsto\left(u^1(z)...u^N(z)    \right)
$
of $\left({\cal U},\Lambda
  \right)$
such that
the initial values
$
\left(
\displaystyle
u^1(z_0).....u^N(z_0),
\U^{\bullet}{}^1(z_0)
.....\U^{\bullet}{}^N(z_0)
  \right)
$
of $\gamma$
yield a metrically
ordinary point of
$\left({\cal U},\Lambda
 \right)$
$u^1$ is not a constant
function
admits the following
first integral:
\begin{eqnarray}
\left(
\U^{\bullet}{}^1
\right)^2
\left( b_1\left(u^1   \right)   \right)
=A_1-\sum_{l=2}^N\frac{A_l}
{a_l\left(u^1   \right)}
\label{eqAa}\\
\left(\U^{\bullet}{}^k   \right)^2
f_k\left(u^k   \right)
\left[a_k\left(u^1   \right)\right]^2
=A_k\quad k=2...N
\label{eqAb}
\end{eqnarray}

Here the $A_k$'s are
suitable complex constants.
\end{lemma}
{\bf Proof:} Let us prove at
first the set of
equations (\ref{eqAb})
corresponding to $k=2...N$.

If $u^k$ is a constant
function, then
$\U^{\bullet}{}^k\equiv 0$
and the k-th equation in
(\ref{eqAb})
holds, with $A_k=0$.

Otherwise,
we could divide the k-th equation in (\ref{equazioninormali}) by $u^k$, this division being lead within the ring of meromorphic functions in a neigbhourhood of $z_0$.
We get
$\displaystyle
2\frac{\U^{\bullet\bullet}{}^k}{\U^{\bullet}{}^k}
+\frac{f^{\prime}_k\left(u^k   \right)}
          {f_k(u^k)}\U^{\bullet}{}^k+
2\frac{a^{\prime}_k\left(u^1   \right)}
          {a_k(u^1)}\U^{\bullet}{}^1=0
$.
Therefore, integrating once,
$
\left(\U^{\bullet}{}^k   \right)^2
f_k\left(u^k   \right)
\left[a_k\left(u^1   \right)\right]^2
=A_k
$
where we have set
$
A_k=
\left(\U^{\bullet}{}^k(z_0)   \right)^2
f_k\left(u^k(z_0)   \right)
\left[a_k\left(u^1(z_0)
 \right)\right]^2
$.
Note that $A_k$ is a well defined complex number, since
$U\left(z_0   \right)=
\left(u^1(z_0)...u^N(z_0)
  \right)$
is a metrically ordinary point.

Let us now prove
(\ref{eqAb}):
we can multiply the first equation of (\ref{equazioninormali}) by $2b_1\left(u^1   \right)\U^{\bullet}{}^1$, since this last function is not everywhere vanishing.

We get
\begin{eqnarray*}
2b_1\left(u^1   \right)
\U^{\bullet}{}^1
\U^{\bullet\bullet}
{}^1+{b_1^{\prime}(u^1)}
\left(\U^{\bullet}{}^1
   \right)^3
-\sum_{l=2}^N {a_l^{\prime}(u^1)f_l(u^l)}
\left(\U^{\bullet}{}^l   \right)^2\U^{\bullet}{}^1=0;
\end{eqnarray*}
by \ref{eqAb} already proved,
$
\left(\U^{\bullet}{}^l
\right)^2=
{A_l}/
[
{f_l(u^l)\left[a_l(u^1)\right
]^2}
]
$,
hence
$$
2b_1\left(u^1
   \right)\U^{\bullet}{}^1
\U^{\bullet\bullet}{}^1+
{b_1^{\prime}(u^1)}
\left(\U^{\bullet}{}^1
   \right)^3
-\sum_{l=2}^N
A_l
\frac
{a_l^{\prime}(u^1)}
{\left[a_l(u^1)\right]^2}
\U^{\bullet}{}^1=0.
$$

Integrating once,
$
b_1\left(u^1   \right)\left(\U^{\bullet}{}^1
   \right)^2+\sum_{l=2}^N
\frac
{A_l}
{a_l\left(u^1   \right)}
=K
$
where
$
K=
b_1\left(u^1(z_0)   \right)\left(\U^{\bullet}{}^1
(z_0)   \right)^2+\sum_{l=2}^N
\frac
{A_l}
{a_l\left(u^1(z_0)   \right)}.
$

Dividing by $b_1\left(u^1   \right)$, keeping into account that $b_1\left(u^1(z_0)   \right)\not=0$ (due to the metrical ordinariness of the initial point of the geodesic) and eventually
setting
$\displaystyle A_1=K/b_1\left(u^1(z_0)   \right)$
ends the proof.
\QUAN
\begin{lemma}
\label{integraleprimo2}
Every element of geodesic
$\displaystyle
z\mapsto\left(u^1...u^N    \right)
$
of $\left({\cal U},\Lambda   \right)$
such that
the initial values
$
\left(
\displaystyle
u^1(z_0).....u^N(z_0),
\U^{\bullet}{}^1(z_0).....\U^{\bullet}{}^N(z_0)
  \right)
$
of $\gamma$
yield a metrically ordinary
point of
$\left({\cal U},\Lambda
\right)$ and $u^1$
 is a constant function
admits the following first integral:
\begin{equation}
\left(\U^{\bullet}{}^k
   \right)^2
f_k\left(u^k   \right)
=A_k\quad k=2...N.
\label{heartsuit}
\end{equation}
Here the $A_k$'s are suitable complex constants.
\end{lemma}
{\bf Proof:}
If $u^k$ is a constant
function, then
$\U^{\bullet}{}^k\equiv 0$
and the k-th equation in
\ref{heartsuit}
holds, with $A_k=0$.

Otherwise, we could
divide the k-th equation in
(\ref{equazioninormali})
by $u^k$, this division
being lead within the ring
of meromorphic functions in
a neigbhourhood of $z_0$.

By keeping into account
that $\U^{\bullet}{}^1\equiv 0$
 we get:
$\displaystyle
2\frac{\U^{\bullet\bullet}{}^k}{\U^{\bullet}{}^k}
+\frac{f^{\prime}_k
\left(u^k   \right)}
          {f_k(u^k)}
\U^{\bullet}{}^k=0
$.
Therefore, integrating once,
$
\left(\U^{\bullet}{}^k   \right)^2
f_k\left(u^k   \right)
=A_k
$,
where we have set
$
A_k=
\left(\U^{\bullet}{}^k(z_0)   \right)^2
f_k\left(u^k(z_0)   \right)
$.
Note that $A_k$ is a well
defined complex number, since
$U\left(z_0   \right)=
\left(u^1(z_0)...u^N(z_0)
   \right)$ is a metrically ordinary point: this fact ends the proof.
\QUAN
\begin{remark}
\label{radiciquadrate}
In the following we shall be concerned with 'extracting square roots' of nonvanishing elements, or germs, of holomorphic functions at some points in the complex plane: more precisely, let $\left(U,\Psi   \right)$ be a never vanishing HFE: then there exist two HFE's
$\left(U,\Xi_1   \right)$ and $\left(U,\Xi_2   \right)$
such that $\Xi_1^2=\Psi$ and $\Xi_2^2=\Psi$ on $U$: the Riemann surfaces of
$\left(U,\Xi_1   \right)$
and $\left(U,\Xi_2   \right)$
are isomorphic, since
{\bf
 either
}
the Riemann surface $\left(R,p,i,\widetilde{U}   \right)$ of
$\left(U,\Psi   \right)$ is such that $\widetilde{U}$ is never vanishing, nor has it got any poles; then the Riemann surfaces of $\left(U,\Xi_1   \right)$, $\left(U,\Xi_2   \right)$  and $\left(U,\Psi   \right)$ are all isomorphic,
{
\bf or}
the Riemann surface $\left(R,p,i,\widetilde{U}   \right)$ of
$\left(U,\Psi   \right)$ is such that there exists some point $p\in R$ such that $\widetilde{U}(p)=0$ or such that $\widetilde{U}$ has a pole in $p$: then the function elements
$\left(U,\Xi_1   \right)$ and $\left(U,\Xi_2   \right)$ are connectible, hence their Riemann surfaces are again isomorphic.
The same argument
 could be applied without changes to the Riemann surfaces of
the HFE's
$\left(U,\int\Xi_1   \right)$ and $\left(U,\int\Xi_2   \right)$.
\end{remark}
\begin{definition}
\label{coercive}
A meromorphic warped product
$$
\displaystyle
{\cal U}={\cal U}_1\times_{a_2(u^1)}
{\cal U}_2\times_{a_3(u^1)}
{\cal U}_3\times
........
\times_{a_N(u^1)}{\cal U}_N
$$
of complex planes or one-dimensional unit balls with metric
$$
\Lambda\left(u^1.....u^N\right)=
b_1(u^1)\,du^i\odot du^i+
\sum_{i=2}^N a_i(u^i)f_i(u^i)
\,du^i\odot du^i,
$$
where $b_1$, the $a_k$'s and the $f_k$'s are
nonzero meromorphic functions
is {\bf coercive} provided that, for every metrically ordinary
 point  $\displaystyle X_0=\left(x_0^1...x_0^N \right)$ and
\begin{itemize}
\item
for every n-tuple
$\left(A_1...A_N   \right)
\in\CI^N$
such that
$
b_1(x_0^1)\not=0$
and
$
A_1-\sum_{l=2}^N\frac
{
A_l}
{
a_l(x_0^1)}\not=0
$
and for each one of the two HFG's
$\displaystyle{\alph_1}$
and
$\displaystyle{\alph_2}$
such that
$$
\left({\alph_i}\right)^2=
\left[
\frac{1}{b_1}
\left(
A_1-\sum_{l=2}^N\frac
{A_l}
{a_l}
\right)
\right]
_{x^1_0}
\quad i=1,2,
$$
the Riemann surface $\displaystyle\left(S_1,\pi_1,j_1,\Phi_1,{\cal U}   \right)$
of both the HFG's
(see remark \ref{radiciquadrate})
\begin{equation}
\left[
\int_{x_0}^{u^1}
{
\frac{d\,\eta}
{\alph_i(\eta)}
}
\right]_{x_0^1}
\quad i=1,2;\label{ciuno}
\end{equation}
is such that $\CI\setminus\Phi_1(S_1)$ is a finite set;
\item
for each $k$, $2\leq k\leq N$
and for each one of the two HFG's
$\displaystyle{\phi_{k1}}$
and
$\displaystyle{\phi_{k2}}$
such that
$$
\left({\phi_{ki}}   \right)^2=
\left[
f_k
\right]_{x_0^1},
\quad i=1,2
$$
the Riemann surface $\displaystyle\left(S_k,\pi_k,j_k,\Phi_k,{\cal U}   \right)$
of both the HFG's
(see remark \ref{radiciquadrate})
\begin{equation}
\displaystyle\left[\int_{x_0^1}^{u^k}
\phi_{ki}(\eta)
\,d\eta\right]_{x_0^1}\ \ i=1,2\label{cidue}
\end{equation}
is such that $\CI\setminus\Phi_k(S_k)$ is a finite set.
\end{itemize}
\end{definition}
\begin{remark}
{Definition \ref{coercive} may be checked
for just one metrically ordinary point $X_0$: this is proved in lemma \ref{isom}; moreover, we may assume,without loss of generality
$X_0=0$:}
were not, we could carry it into $0$ by
applying an automorphism of
${\cal U}$, that is to say a direct product of automorphisms of the unit ball or of the complex plane, according to the nature of each ${\cal U}_i$.
Then a simple pullback procedure would yield back the initial situation: {in the following we shall understand this choice}.
\end{remark}
In the following lemma we shall use the 'square root' symbol in the meaning of definition \ref{coercive}, or remark \ref{radiciquadrate}: in other words, given a HFG, which is not vanishing at some point, it should denote any one of the two HFG's yielding it back when squared.
\begin{lemma}
\label{isom} For every metrically ordinary point
$\displaystyle\left (\xi^1...\xi^N   \right)$ of ${\cal U}$ and
every n-tuple $\displaystyle\left(A_1...A_N   \right)\in\CI^N$
such that $ b_1(x_0^1)\not=0$, $ A_1-\sum_{l=2}^N\frac { A_l} {
a_l(x_0^1)}\not=0$, $ \ b_1(\xi^1)\not=0$ and $
A_1-\sum_{l=2}^N\frac { A_l} { a_l(\xi^1)}\not=0$, set
$\Psi(\eta):= \displaystyle{ A_1-\sum_{l=2} ^N\frac {A_l}
{a_l(\eta)}}$: then the Riemann surfaces of the HFG's
$\displaystyle \int_{\xi_1}^{u^1} \sqrt{\displaystyle b_1
(\eta)/\Psi(\eta)} \,d\eta \quad $ at ${\xi_1}$ and $
\int_{0}^{u^1} \sqrt{\displaystyle b_1(\eta)/\Psi(\eta)} \,d\eta
\quad $ at ${0} $ are isomorphic: moreover so are, for each $k$,
those of $\displaystyle \int_{\xi_k}^{u^k} \sqrt{ f_k(\eta)
\,d\eta } $ at ${\xi_k} $ and $\displaystyle \int_{0}^{u^k} \sqrt{
f_k(\eta) \,d\eta } $ at ${0} $.
\end{lemma}
{\bf Proof:} The statement easily
follows from the fact that
those germs are connectible.
\QUAN
Here is the main result
concerning warped
products of Riemann surfaces:
\begin{theorem}
\label{teoremaprincipale}
A meromorphic warped product
$$
\displaystyle
{\cal U}={\cal U}_1\times_{a_2(u^1)}
{\cal U}_2\times_{a_3(u^1)}
{\cal U}_3\times
........
\times_{a_N(u^1)}{\cal U}_N
$$
of complex planes or one-dimensional unit balls with metric
$$
\Lambda\left(u^1.....u^N\right)=
b_1(u^1)\,du^1\odot du^1+
\sum_{i=2}^N a_i(u^1)f_i(u^i)
\,du^i\odot du^i,
$$
is geodesically complete if and only if it is coercive.
\end{theorem}
{\bf Proof:}
a) Suppose that $\displaystyle{\cal U}$ is {\TTT coercive} and that $U$, defined by
$
z\mapsto\left(u^1...u^N    \right),
$
is an element of geodesic, defined in a neighbourhood of $0$ in the complex plane and such that $\left(u^1(0)...u^N(0)    \right)$ is a metrically ordinary point; moreover, let
$
\left(\U^{\bullet}{}^1(0)...\U^{\bullet}{}^N(0)    \right)
$
be the initial velocity of $U$.

Suppose at first that $\displaystyle z\mapsto u^1$ is a constant function (hence $\U^{\bullet}{}^1(0)=0\,$): then, by lemma \ref{integraleprimo2}, the equations of $U$ are
\begin{equation}
\label{riportointpr2}
\left(\U^{\bullet}{}^k
   \right)^2
f_k\left(u^k   \right)
=A_k\quad\  k=2...N,
\end{equation}
where the $A_k$'s are
suitable complex constants;
here $u^1\equiv A_1$.

Now the Riemann surface of the HFE $\displaystyle z\mapsto {u}^1$ is trivially isomorphic to $\displaystyle \left(\CI,\id,\id,A_1 \right)$; if  $\displaystyle A_k=0$ the Riemann surface of $\displaystyle z\mapsto {u}^k$ is isomorphic to $\displaystyle \left(\CI,\id,\id,A\right)$ for some complex constant $A$; if $\displaystyle A_k\not=0$ we could rewrite the k-th equation of (\ref{riportointpr2}) in the form:
\begin{equation}
\label{riscritta}
\frac{1}{B_k}\int_{u^k(0)}^{u^{k}} \phi(\eta)\,d\eta
=\,z,
\end{equation}
where $\phi_k^2=f_k$ and
$B_k^2=A_k$, the choice
of $\phi_k$ and $B_k$
being made in such a
way that
$\displaystyle
\U^{\bullet}{}^k(0)=
\frac{B_k}{\phi_k(0)}$.

By hypothesis, the Riemann surface $\displaystyle\left(S_k,\pi_k,j_k,\Phi_k   \right)$
of the HFG $\displaystyle
\int_0^{u^k}\phi_k\,d\eta$ at ${0}$
is such that $\displaystyle\CI\setminus\Phi_1(S_1)$ is a finite set; by lemma \ref{isom} the
Riemann surface of the HFG
$\displaystyle
\int_{u^k(0)}^{u^k}\phi_k\,d\eta$
at
$
{u^k(0)}$
is isomorphic to $\displaystyle\left(S_k,\pi_k,j_k,\Phi_k   \right)$; but, by (\ref{riscritta}), the germs
$\displaystyle \u^k_{z=0}$ and
$\displaystyle
\int_{u^k(0)}^{u^k}\phi_k\,d\eta$
at
${u^k(0)}$
are each one inverse of the other; hence, by lemma \ref{inverse} the Riemann surface of $\displaystyle \u^k_{z=0}$ is complete; this eventually implies that the Riemann surface of the element
$
z\mapsto\left(u^1...u^N    \right)
$
is complete too: this fact ends the proof of a) in the case that $u^1$ is a constant function.

On the other side, suppose that $u^1$ is not a constant function: then, by lemma \ref{integraleprimo1}, the equations of $U$ are
\begin{equation}
\label{riportointpr1}
\cases{
\left(\U^{\bullet}{}^1 \right)^2
\left( b_1\left(u^1   \right)   \right)
=A_1-\sum_{l=2}^N\frac{A_l}
{a_l\left(u^1   \right)}
\quad\cr
\left(\U^{\bullet}{}^k   \right)^2
f_k\left(u^k   \right)
\left[a_k\left(u^1   \right)\right]^2
=A_k\quad k=2...N.
}
\end{equation}
for suitable complex
constants $A_1...A_N.$

Consider now the germ $\displaystyle z\mapsto u^1$ in $z=0$:
rewrite the first equation of (\ref{riportointpr1}) in the form:
\begin{equation}
\label{riscritta2}
\int_{u^1(0)}^{u^{1}}
\frac{\displaystyle d\eta}
 {\displaystyle\alph(\eta)_{u^1(0)}}
=\,z,
\end{equation}
where
$\displaystyle
\left(
\alph(\eta)_{u^1(0)}
\right)^2=
\displaystyle
\left[{A_1-\sum_{l=2}^N
A_l
/
a_l(\eta)}
    \right]/
{b_1(\eta)}
$
in a neighbourhood of $z=0$,
the choice of the square root  $\alph_k$ being
made in such a way that
$
\displaystyle \alph_{u^1(0)}\left(u^1(0)\right)=1/\U^{\bullet}{}^1(0)
$.

Denote now  by $\alph_{u=0}$ the HFG
defined by setting
$$
\left(
\alph_{0}
\right)^2=
\left[
\frac
{1}
{\displaystyle b_1}
\left(
{\displaystyle
A_1-\sum_{l=2}^N\displaystyle\frac
{A_l}
{a_l}}
\right)
\right]_{0}
,
$$
the choice of the 'square root' $\alph_0$ being arbitrary.

By hypothesis, the Riemann surface $\displaystyle\left(S_1,\pi_1,j_1,\Phi_1   \right)$
of the HFG $\displaystyle
\int_0^{u^1}1/\alph_0$
at $0$
is such that $\displaystyle\CI\setminus\Phi_1(S_1)$ is a finite set.

By lemma \ref{isom} the Riemann surfaces of $\displaystyle
\int_0^{u^1}1/\alph_0$ (at $0$) and of $\displaystyle\int_{u^1_0}^{u^1}1/\alph_0
$ (at ${u^1_0}$)
are both isomorphic to $\displaystyle\left(S_1,\pi_1,j_1,\Phi_1   \right)$; but, by (\ref{riscritta}), the germs
$\displaystyle \u^1_{z=0}$ and
$\displaystyle\int_0^{u^1}1/\alph_0$
(at $
{u^1(0)}$) are each one inverse of the other;
hence, by lemma \ref{inverse} the Riemann surface of $\displaystyle \u^1_{z=0}$ is complete.

Let now $\displaystyle 2\leq k\leq N$: if $\displaystyle A_k=0$ the Riemann surface of $\displaystyle z\mapsto {u}^k$ is isomorphic to $\displaystyle \left(\CI,\id,\id,A\right)$ for some complex constant $A$; if $\displaystyle A_k\not=0$ we could rewrite the k-th equation of (\ref{riportointpr1}) in the form:
\begin{equation}
\label{riscritta3}
\int_{u^k(0)}^{u^{k}} \phi(\eta)\,d\eta
=\,
\int_0^z\frac{B_k\,dz}{a_k\left(u^1   \right)},
\end{equation}
where $\phi_k^2=f_k$
and $B_k^2=A_k$,
the choice of $\phi_k$ and
$B_k$ being made in such
a way that
$\displaystyle \U^{\bullet}
{}^k(0)\,\phi\left(
u^k(0)   \right)\,
a_k\left(u^1   \right)=
{B_k}$.

Denote now  by $\displaystyle[\varphi_k]_{u^k=0}$ the HFG
defined by setting
$\displaystyle[\varphi_k]_{u^k=0}^2=
\left[
f_k
\right]_{u^k=0}
$,
the choice of the "square root" $\displaystyle[\varphi_k]_{u^k=0}$ being arbitrary.

By hypothesis, the Riemann surface $\displaystyle\left(S_k,\pi_k,j_k,\Phi_k   \right)$
of the HFG $\displaystyle
\int_0^{u^k}\varphi_k$ (at ${0}$)
is such that $\displaystyle\CI\setminus\Phi_1(S_1)$ is a finite set; moreover, by lemma \ref{isom} the Riemann surfaces of the HFG $\displaystyle
\int_{u^k(0)}^{u^k}\phi_k\,d\eta$
( at ${u^k(0)}$)
is isomorphic to
$\displaystyle\left
(S_k,\pi_k,j_k,\Phi_k
   \right)$; but, by
(\ref{riscritta3}) the germs
$\displaystyle
\left[z\rightarrow
\u^k\right]_{z=0}$,
$\displaystyle
\int_{u^k(0)}^{u^k}
\phi_k\,d\eta$ (at ${u^k(0)}$
and
$\displaystyle
z\rightarrow
\int_0^z\frac{B_k}
{a_k\left(u^1(\zeta)
  \right)}
\,d\zeta$
( at ${z=0}$)
satisfy, in the above order, the hypotheses of lemma \ref{quasiinverse};
moreover, the Riemann surface with logarithmic singularities of $\displaystyle
\int_{u^k(0)}^{u^k}\phi_k\,d\eta
$ ( at ${u^k(0)}$) is complete, since the one of
$\displaystyle\left[\phi_k\right]_{u^k(0)}$ is complete without logarithmic singularities.

Therefore the Riemann surface
 with logarithmic
singularities of
$\displaystyle \u^k_{z=0}$
is complete;
this eventually implies
that the Riemann surface
with logarithmic singularities
of the element
$
z\mapsto\left(u^1...u^N    \right),
$
is complete too:
this fact ends the proof of a).
\vskip0,2truecm
Vice versa, suppose that
$
\displaystyle
{\cal U}={\cal U}_1\times_{a_2(u^1)}
{\cal U}_2\times_{a_3(u^1)}
{\cal U}_3\times
........
\times_{a_N(u^1)}{\cal U}_N
$
is not {\TTT coercive}: then
{\bf either}
there exists a complex
n-tuple $\displaystyle \left(A_1...A_N   \right)\in\CI^N$ such that
$
\displaystyle
b_1(x_0^1)\not=0$,
$
A_1-\sum_{l=2}^N\frac
{\displaystyle
A_l}
{\displaystyle
a_l(x_0^1)}\not=0
$
and for each one of the two HFG's
$\displaystyle{\alph_1}$
and
$\displaystyle{\alph_2}$
such that
$$
\left({\alph_i}\right)^2=
\left[
\frac{1}{b_1}
\left(
A_1-\sum_{l=2}^N\frac
{A_l}
{a_l}
\right)
\right]
_{0}
\quad i=1,2,
$$
the Riemann surface $\displaystyle\left(S_1,\pi_1,j_1,\Phi_1   \right)$
of both the HFG's
(see remark \ref{radiciquadrate})
$
\int_{x_0}^{u^1}
{
\frac{d\,\eta}
{\alph_i(\eta)}
}
$
( at $
{x_0^1}
(\  i=1,2)
$)
is such that $\CI\setminus\Phi_1(S_1)$ is
an infinite set;
{\bf or}
there exists $k$,
$2\leq k\leq N$
such that, for each one of
the two HFG's
$\displaystyle
\left[\phi_{k1}\right]_0$
and
$\displaystyle
\left[\phi_{k2}\right]_0$
such that
$
\displaystyle\left[\phi_{ki}\right]_0=
\left[
f_k
\right]_0,
\ (i=(1,2))
$
the Riemann surface $\displaystyle\left(S_k,\pi_k,j_k,\Phi_k   \right)$
of both the HFG's
(see remark \ref{radiciquadrate})
$\displaystyle\left[\int_0^{u^k}
\phi_ki(\eta)
\,d\eta\right]_{0}\ \ i=1,2$
is such that $\CI\setminus\Phi_1(S_1)$ is an infinite set.

In the first case the
geodesic element
$
z\mapsto U=\left(u^1...u^N    \right)
$
starting from $0$ with velocity
$\displaystyle\left(L_1...L_N    \right)$, such that
$$\displaystyle
L_1^2=
\frac
{1}
{b_1(0)}
\left(
A_1-\sum_{l=2}^N\frac
{A_l}
{a_l(0)}
\right)
,\ L_k^2=\frac{A_k}{f_k(0)a_k(0)},\
k=2...N
,$$
satisfies the equation
$
\int_0^{u^1}
\frac
{\displaystyle d\eta}
{\displaystyle \alph_i(\eta)}
=z
$,
$i=1,2$;
by lemma \ref{inverse}, this
fact implies that $\left[
\displaystyle z\mapsto
u^1
\right]_{0}  $
has an incomplete Riemann
surface,
hence  the same holds
about $\displaystyle z\mapsto
U$ too.

Consider now the second case: first construct a geodesic element
$
z\mapsto U=\left(0...u^k...0    \right),
$
with all components
which have to be constant
functions except
$\displaystyle u^k, k\geq 2$
(this element is easily
seen to exist).

Now recall lemma
\ref{integraleprimo2}
to conclude that
$
\displaystyle
z\mapsto u^k
$
satisfies,
in a neighbourhood
of $z=0$ the equation
$
\displaystyle
\frac
{1}
{C_k}
{\int_0^{u^k} \phi_{ki}(\eta)\,d\eta}=z,
$
for a suitable complex constant $A_k$; therefore its Riemann surface is incomplete by lemma \ref{inverse}; this fact ends the proof.
\QUAN
\begin{definition}
\label{directbiholom}
Let ${\cal U}$ and ${\cal V}$ be meromorphic warped products of complex planes and unit balls;  ${\cal U}$ and ${\cal V}$ are {\bf directly biholomorphic} provided that they are biholomorphic under a direct product of biholomorphic functions between each  ${\cal U}_i$ and each ${\cal V}_i$.
\end{definition}
\begin{remark}
\label{biholom}
Definition \ref{coercive} is invariant under direct biholomorphism (see definition \ref{directbiholom}): in other words, if ${\cal U}$ and ${\cal V}$ are directly biholomorphic, then
${\cal U}$ is coercive if and only ${\cal V}$ is too: this is a simple consequence of 'changing variable' in integrals \ref{ciuno} and \ref{cidue}.
\end{remark}
Therefore, we could yield the following
\begin{definition}
\label{equicoercive}
An equivalence class $\left[{\cal U}    \right]$ of meromorphic warped products of complex planes and unit balls, consisting of mutually directly (see definition \ref{directbiholom}
) biholomorphic elements is {\bf coercive}
provided that any one of its representatives is coercive.
\end{definition}

Our goal is now to extend definitions \ref{coercive} and \ref{equicoercive} to warped products containg some $\PI^1$'s among their factors.

Keeping into account remark \ref{biholom}, this could be readiliy pursued: indeed, consider a meromorphic warped product
$$
\displaystyle
{\cal U}={\cal U}_1\times_{a_2(u^1)}
{\cal U}_2\times_{a_3(u^1)}
{\cal U}_3\times
........
\times_{a_N(u^1)}{\cal U}_N
$$
of Riemann spheres, complex planes or one-dimensional unit balls with metric
$$
\Lambda\left(u^1.....u^N\right)=
b_1(u^1)\,du^i\odot du^i+
\sum_{i=2}^N a_i(u^i)f_i(u^i)
\,du^i\odot du^i.
$$

Let $\displaystyle L\subset \{1...N\}$ be the set of indices such that $\displaystyle {\cal U}_l\simeq \PI^1$ for each $l\in L$.

\begin{definition}
Let $\displaystyle
Y=\left(y^1...y^N \right)
\in{\cal U}$:
then $\left(Y,L    \right)$
is a {\bf principal multipole}
 of ${\cal U}$ provided that
$
b_1(y^1)=\infty$ and
$
f_l(y^l)=\infty$
  for each
$l\in L\setminus \{1\}$.
\label{multipole}
\end{definition}
\begin{definition}
A meromorphic warped product
$$
\displaystyle
{\cal U}={\cal U}_1\times_{a_2(u^1)}
{\cal U}_2\times_{a_3(u^1)}
{\cal U}_3\times
........
\times_{a_N(u^1)}{\cal U}_N
$$
of Riemann spheres, complex planes or one-dimensional unit balls with metric is {\bf partially projective} if some one of its factors is biholomorphic to the Riemann sphere $\PI^1$.
\end{definition}
\begin{definition}
\label{multicoer}
A partially projective
warped product
$\displaystyle {\cal U}=
\prod_{i=1}^N{\cal U}_i$
is {\bf coercive in
opposition to the principal
multipole }
$\left(Y,L    \right)$ if,
set
$
{\cal W}_i=
{\cal U}_i$  if  $i\not\in L$,
$
{\cal W}_i={\cal U}_i\setminus \{y^i\}$
$ if i\in L
$,
then $\displaystyle \prod_{i=1}^N {\cal W}_i$ is coercive in the sense of definition \ref{equicoercive}, that is to say, belongs to a coercive equivalence class with respect to direct biholomorphicity.
\end{definition}

\subsection{Warped product of Riemann surfaces}

Consider now a warped product of Riemann surfaces
$$
\displaystyle
{\cal S}={\cal S}_1\times_{a_2}
{\cal S}_2\times_{a_3}
{\cal S}_3\times
........
\times_{a_N}{\cal S}_N,
$$
where each ${\cal S}_i$ is endowed with meromorphic metric $\lambda_i$: ${\cal S}$'s metric $\Lambda$ is defined by setting
$$
\Lambda=\lambda_1+\sum_{k=2}^N a_k\lambda_k,
$$
where each $a_k$ is a meromorphic function on ${\cal S}_i$.

\begin{theorem}
\label{rivestimento}
${\cal S}$ admits universal covering $\displaystyle \Psi : {\cal U}\rightarrow {\cal S}$, where ${\cal U}$ is a direct product
of Riemann spheres, complex planes or one-dimensional unit balls: this universal covering is unique up to direct biholomorphisms.
\end{theorem}
{\bf Proof:} This is a simple consequence of Riemann's uniformization theorem.
\QUAN
Now ${\cal U}$ could be endowed with the pull-back meromorphic metric $\displaystyle \Psi^*\Lambda$, hence ${\cal U}$ itself results in a meromorphic warped product.
\begin{definition}

\label{rieco}
The manifold ${\cal S}$ is
{\bf totally unelliptic}
provided that none of
the ${\cal S}_i$ is elliptic;
{\bf $L$-elliptic}
provided that there exists
a nonempty set of indices
$L$ such that ${\cal S}_l$
is elliptic if and only
if $l\in L$.
\end{definition}
\begin{definition}
Let ${\cal S}$ be a $L$-elliptic warped product, with universal covering $\Psi : {\cal U}\rightarrow{\cal S}$: then $\left(Z,L    \right)$ is a principal multipole for ${\cal S}$ provided that $Z\in {\cal S}$ and each $Y\in\Psi^{-1}\left(Z    \right)$ is a principal multipole for ${\cal U}$.
\end{definition}
\begin{definition}
A totally unelliptic
warped product of Riemann
sur\-faces is {\bf coercive}
provided that its universal
 covering is coercive in
the sense of definition
\ref{equicoercive}.
A $L$-elliptic warped
product of Riemann surfaces
is {\bf coercive in opposition
 to the principal multipole}
$\displaystyle\left(Z,L
  \right)$ provided that
its universal covering
${\cal U}$ is coercive
in opposition to each
principal multipole
$\displaystyle\left(Y,L
  \right)$ as $Y$ runs
over $\Psi^{-1}(Z)$.
\end{definition}
\begin{theorem}
\label{unell}
A totally unelliptic warped product of Riemann surfaces ${\cal S}$ is geodesically complete if and only if it is coercive.
\end{theorem}
{\bf Proof:} Let $\Psi : {\cal U}\rightarrow
{\cal S}$ be the universal covering of ${\cal S}$:
by definition \ref{rieco} ${\cal U}$ is coercive, hence geodesically complete by theorem \ref{teoremaprincipale}.

Let now $\gam$ be a germ of geodesic in ${\cal S}$, starting at a metrically ordinary point: since $\Psi$ is a local isometry, there exists a germ $\bet$ of geodesic in ${\cal U}$, starting at a metrically ordinary point, such that $\displaystyle\gam=\Psi\circ\bet$.

By definition of completeness, the Riemann surface with logarithmic singularities $\displaystyle\left(\Sigma,\pi,j,B,{\cal U}    \right)$ of $\bet$ is such that $\CI\setminus \pi\left(\Sigma    \right)$ is a finite set; moreover, $\displaystyle\left(\Sigma,\pi,j,\Psi\circ B,{\cal S}    \right)$ is an analytical continuation, with logarithmic singularities, of $\gam$.

This proves that, if $\displaystyle\left(\widetilde \Sigma,\widetilde \pi,\widetilde{\hbox{\j}},G,{\cal S}    \right)$ is the Riemann surface with logarithmic singularities of $\gam$, then $\PI^1\setminus \widetilde \pi\left(\widetilde \Sigma    \right)$ is a finite set too, hence ${\cal S}$ is geodesically complete.

On the other side, if ${\cal S}$ admits an incomplete germ of geodesic $\gam$, starting at a metrically ordinary point, then there exists  an incomplete germ of geodesic $\bet$ in ${\cal U}$, starting at a metrically ordinary point, such that
$\displaystyle\gam=\Psi\circ\bet$; this means by theorem \ref{teoremaprincipale},
that ${\cal U}$ is not coercive; eventually, by definition \ref{rieco}, ${\cal S}$ is not coercive: this fact ends the proof.
\QUAN
\begin{theorem}
A $L$-elliptic warped product of Riemann surfaces ${\cal S}$ is geodesically complete if and only if) it is coercive in opposition to some principal  multipole.
\end{theorem}
{\bf Proof:} Suppose that ${\cal S}$ is coercive in opposition to some principal  multipole $\displaystyle\left(Z,L    \right)$: then, by theorem \ref{unell},  ${\cal S}$ is coercive in opposition to $\displaystyle\left(Z,L    \right)$ if and only if ${\cal S}\setminus Z$ is geodesically complete;
since $Z$ is not metrically ordinary, ${\cal S}$ is geodesically complete.

On the other hand, suppose that ${\cal S}$ admits an incomplete geodesic
$\displaystyle\left(\Sigma,\pi,j,\gamma,{\cal S}    \right)$: let $\left(Z,L    \right)$ be a principal multipole of ${\cal S}$ wich is known to exist; set $R=\gamma^{-1}\left({\cal S}\setminus Z    \right)\subset\Sigma$.

Now $\displaystyle\left(R,\pi\vert_R,j,\gamma\vert_R,{\cal S}\setminus Z    \right)$ is an incomplete geodesic of ${\cal S}\setminus Z$: this fact implies that ${\cal S}\setminus Z$ is not geodesically complete, hence it is not coercive, that is to say, ${\cal S}$ is not coercive in opposition to $\left(Z,L    \right)$.

The arbitrariness of $Z$ allows us to conclude the proof.
\QUAN

\subsubsection{Examples}
We shall now show a wide class of warped products sharing all characteristics defining coercivity: they will hence result in being geodesically complete.

We recall, without proof, the following results from the theory of meromorphic functions (see \cite{nevanlinna} or \cite{hayman}):
\begin{theorem}
\label{mero}
A meromorphic function in the complex plane takes all $\PI^1$'s values but at most two ones;
a meromorphic function  in the unit disc, whose characteristic function $T$ is such that the ratio
$
T(r)/log(1-r)
$
is unlimited as $r\rightarrow 1$, takes all $\PI^1$'s values but at most two ones.
\end{theorem}

In the following we shall need some technicalities
 from integral calculus,
hence we state:
\begin{proposition}
\label{integrale}
If $\Delta:=b^2-4ac$
then
$
\left[
\int \frac{\displaystyle d\,\eta}
{\sqrt{a\eta^2+b\eta+c}}
\right]_0
$ equals one  and only
one (up to additive constants)
of the following expressions, in a neighbourhood of $0$:
\begin{eqnarray*}
&\hbox{\bf (A)}&
\cases{
\displaystyle
\left[
\frac{1}
{\sqrt{a}}\log
\left(\eta+\frac{b}{2a}+
\sqrt{\eta^2+\frac{b}{a}
\eta+\frac{c}{a}}
    \right)
\right]_0
\cr
\hbox{\sl
the same branch of $\sqrt{}$,
any branch of log}\cr
\hbox{\sl if $a\not=0$
and $\Delta\not=0$}
}
\ \hbox{\bf (B)}
\cases{
\left[ \frac{1}{\sqrt{a}}
 \log\left(\eta+\frac{b}{2a}
\right)\right]_0 \cr \hbox{\sl any branch of log}\cr \hbox{\sl if
$a\not=0$ and $\Delta=0$}\cr
}\\
&\hbox{\bf (C)}&
%\hbox{\bf (C)}
\cases{
\left[ \frac{\displaystyle 2}{\displaystyle b}\sqrt{b\eta+c}
\right]_0 \cr \hbox{\sl the same branch of $\sqrt{}$} \cr
\hbox{\sl if $a=0$ and $b\not=0$}
}
\hbox{\bf (D)}
\cases{
\left[ \eta/\sqrt{c}\right]_0\cr \hbox{\sl the same branch of
$\sqrt{}$} \cr \hbox{\sl if $a=b=0$.}
}
\end{eqnarray*}
\QUAN
\end{proposition}
Let now $h,f_2...f_N$ be
meromorphic functions
on $\CI$ and $P_2...P_N$
polynomials of degree at
most two.
Consider on $\CI^N$ the meromorphic metric
$$
\Lambda\left(u^1...u^N\right)=
\left[h^{\prime}(u^1)\right]^2 du^1\odot du^1+
\sum_{k=2}^N
\frac{\displaystyle\left[f_k(u^k)\right]^2}
{P_k\left(h(u^1)\right)}\, du^k\odot du^k.
$$
\begin{theorem}
$\displaystyle\left(\CI^N,\Lambda\right)$ is coercive (hence geodesically complete).
\label{esempio}
\end{theorem}
{\bf Proof:}
For every n-tuple
$\displaystyle
 \left(A_1...A_N
 \right)\in\CI^N$ such that
$
h^{\prime}(0)\not=0
$ and
$
A_1-\sum_{l=2}^N
{A_l}
{P_l(0)}\not=0
$, set
$\Psi(x)=A_1-\sum_{l=2}^N
{A_l}
{P_l(x)} $.
There holds
$
\displaystyle\int_0^{u^1}
(\Psi\circ h(\eta))^{-1/2}
h^{\prime}(\eta)d\,\eta=
\Phi\left(h(u^1)    \right)
$,
where $\Phi$ is one (depending on the constants $A_1...A_N$) of the HFG's on the right hand member of proposition \ref{integrale}.
This fact shows that the maximal analytical continuation of
$u^1\rightarrow\Phi\left(h(u^1)    \right)$ takes all $\PI^1$'s values but a
finite number, because so does the meromorphic function $h$ (see theorem \ref{mero}).

Moreover,
for each $k$, $2\leq k\leq N$,
each one of the two HFG's
$
\pm
\left[
f_k
\right]_0,
$
could be continuated to $\pm f_k$ which, by theorem \ref{mero}, takes all values but at most two ones.
\QUAN
\begin{remark}
 Extending the validity of preceding example to the partially projective case is straightforward.
\end{remark}

Let now $S_i,\ i=1..N$ be Riemann surfaces, which we suppose for simplicity to be
parabolic or hyperbolic, $p_i\colon{\cal U}_i\rightarrow S_i$ their universal covering,
where each ${\cal U}_i$ is isomorphic either to the unit disc or to the complex plane;
finally, let $\phi_i$ be meromorphic functions such that $\phi_1\circ p_1$ and
$(\phi_i\circ p_i)^{\prime},\ i=1..N$ take all complex values but at most a finite number (the hypothesis on $\phi_i\circ p_i$ could be weakened; even dropped, if $S_i$ is parabolic: see \cite{hayman}).

Moreover, let $a_i,\ b_i,\ c_i,\ i=1..N $ be complex numbers such that, for each $i$, $a_i\not=0$ or  $b_i\not=0$ or $\ c_i\not=0$.

Set
$
S=\prod_{i=1}^N S_i$,
$
{\cal U}=\prod_{i=1}^N=
{\cal U}_i$
and
$
p=(p_1....p_N)$;
consider the meromorphic metric
$$
\Lambda=d\phi_1\odot d\phi_1+\sum
_{i=1}^N \frac {d\phi_i\odot d\phi_i}{a_i\phi_1^2+b_i\phi_1+c_i}.
$$
\begin{theorem}
$({\cal U},\Lambda)$ is coercive (hence geodesically complete).
\end{theorem}
{\bf Proof:} By pulling back $\Lambda$ with respect to the universal covering $p$ we get
$$
p^*\Lambda(z^1...z^N)=
\left[(\phi_1\circ p_1)^{\prime}\right]^2dz^1\odot dz^1\ +\sum
_{i=1}^N \frac
{\left[(\phi_i\circ p_i)^{\prime}\right]^2dz^i\odot dz^i}
{a_i(\phi_1\circ p_1)^2+b_i\phi_1\circ p_1+c_i}.
$$
We claim that $\left({\cal U},p^*\Lambda\right)$ is coercive:
in fact, for every n-tuple $\displaystyle \left(A_1...A_N   \right)\in\CI^N$ such that
$$
\cases{
(\phi_1\circ p_1)^{\prime}(0)\not=0\cr
A_1-\sum_{l=2}^N
{A_l}
{a_i(\phi_1\circ p_1)^2+
b_i\phi_1\circ p_1+c_i}
\Big\vert_0\not=0,
}
$$
set
$\Psi(x
):=A_1-\sum_{l=2}^N
{A_l}
{a_i(x)^2+
b_i x+c_i}
$,
there holds
$$
\int_0^{u^1}
(\Psi(\phi_1\circ p_1(\eta)))
^{-1/2}
(\phi_1\circ p_1)^
{\prime}(\eta)d\,\eta
=
\int_{\phi_1\circ p_1(0)}^{\phi_1\circ p_1(u^1)}
=
\Phi\left(\phi_1\circ p_1\right),
$$
where $\Phi$ is one (depending on the constants $A_1...A_N$) of the holomorphic function
germs on the right hand member of proposition \ref{integrale}.

This fact shows that the maximal
analytical continuation
of $u^1\rightarrow\Phi
\left(\phi_1\circ p_1(u^1)
    \right)$ takes all $\PI^1$'s values but a finite number, because so
does the meromorphic function
 $\phi_1$ and hence
$\phi_1\circ p_1$;
moreover,
for each $i$, $2\leq i\leq N$,
each one of the two HFG's
$
\pm
\left[(\phi_i\circ p_i)^{\prime}\right]
$
could be continuated to $\pm \left[(\phi_i\circ p_i)^{\prime}\right]$ which, by
assumption, takes all values but at most two ones.
\QUAN

The preceding examples may be readily extended to the following two (alternative) cases, mostly following the outline of the above reasoning:
\def\DI {\sdopp {\hbox{D}}}
\begin{itemize}
\item $\DI^N$ taking place of $\CI^N$ and $h,f_2...f_N$ meromorphic functions on $\DI$ satisfying theorem \ref{mero};
\item $P_2...P_N$ polynomials of degree at most four: similar conclusions may be drawn by means of elliptic integrals.
\end{itemize}
\subsection{Pseudo-Riemannian
warped products}
\begin{definition}
\label{completezzareale}
A pseudo-Riemannian manifold
is {\TTT geodesically complete} provided that it admits a complexification $\M$
such that the Riemann surface,
with logarithmic singularities,
of each (complexified)
geodesic
germ is real-complete.
\end{definition}
\begin{definition}
\label{realcoercive}
A warped product
$$
\displaystyle
{\cal U}={\cal U}_1\times_{a_2(u^1)}
{\cal U}_2\times_{a_3(u^1)}
{\cal U}_3\times
........
\times_{a_N(u^1)}{\cal U}_N
$$
of real intervals, real lines or $\ESSE^1$'s with nondegenerating real-analytic pseudo-Riemannian metric
$$
\Lambda\left(u^1.....u^N\right)=
b_1(u^1)\,du^i\odot du^i+
\sum_{i=2}^N a_i(u^i)f_i(u^i)
\,du^i\odot du^i,
$$
of arbitary signature
is {\bf coercive} provided that, called ${\cal K}$ the canonical complexification $\ERRE^N\rightarrow\CI^N$,
for one (hence every) point $X_0=(x_0^1...x_0^N)$
there holds:
\begin{itemize}
\item
for every n-tuple $\displaystyle \left(A_1...A_N   \right)\in\ERRE^N$ such that
$
\displaystyle
b_1(x_0^1)\not=0$
and
$
A_1-\sum_{l=2}^N\frac
{\displaystyle
A_l}
{\displaystyle
a_l(x_0^1)}\not=0
$
and for each one of the two HFG's
$\displaystyle{\alph_1}$
and
$\displaystyle{\alph_2}$
such that
$$
\left({\alph_i}\right)^2=
{\cal K}\circ
\left[
\frac{1}{b_1}
\left(
A_1-\sum_{l=2}^N\frac
{A_l}
{a_l}
\right)
\right]
_{x_0^1}
\quad i=1,2,
$$
the Riemann surface $\displaystyle\left(S_1,\pi_1,j_1,\Phi_1\right)$
of both the HFG's
(see remark
\ref{radiciquadrate})
\begin{equation}
\left[
\int_{x_0}^{u^1}
{
\frac{d\,\eta}
{\alph_i(\eta)}
}
\right]_{x_0^1}
\quad i=1,2;
\end{equation}
is such that $\ERRE\setminus\Phi_1(S_1)$ is a finite set;
\item
for each $k$, $2\leq k\leq N$
and for each one of the two HFG's
$\displaystyle{\phi_{k1}}$
and
$\displaystyle{\phi_{k2}}$
such that
$$
\left({\phi_{ki}}   \right)^2=
{\cal K}\circ
\left[
f_k
\right]_{x_0^1},
\quad i=1,2
$$
the Riemann surface $\displaystyle\left(S_k,\pi_k,j_k,\Phi_k   \right)$
of both the holomorphic function germs
(see remark \ref{radiciquadrate})
\begin{equation}
\displaystyle\left[\int_{x_0^1}^{u^k}
\phi_{ki}(\eta)
\,d\eta\right]_{x_0^1}\ \ i=1,2
\end{equation}
is such that $\ERRE\setminus\Phi_k(S_k)$ is a finite set.
\end{itemize}
\end{definition}
We confine ourselves in stating the real analogue of our main theorem (the proof is almost identical):
\begin{theorem}
A warped product
$$
\displaystyle
{\cal U}={\cal U}_1\times_{a_2(u^1)}
{\cal U}_2\times_{a_3(u^1)}
{\cal U}_3\times
........
\times_{a_N(u^1)}{\cal U}_N
$$
of real intervals, real lines or $\ESSE^1$'s with nondegenerating real-analytic pseudo-Riemannian metric
$$
\Lambda\left(u^1.....u^N\right)=
b_1(u^1)\,du^i\odot du^i+
\sum_{i=2}^N a_i(u^i)f_i(u^i)
\,du^i\odot du^i,
$$
of arbitary signature is geodesically complete if and only if it is coercive.
\end{theorem}

\subsection{The Clifton-Pohl torus}
Consider now
$\N:=\ERRE^2
\setminus \{0\}$, with the Lorentz metric
${du\odot dv}/({u^2+v^2})
$;
the group $D$ generated by scalar
multiplication by $2$ is a group of isometries of
$\N$;
its action is properly dicontinuous, hence
$\T=\N/D$ is a Lorentz surface.
Topologically, $\T$
is  the closed annulus $1\leq\varrho\leq 2$, with boundaries identified by the action of $D$, i.e. a torus; notwithstanding, $\T$ is geodesically incomplete, since
$t\mapsto\left(1/(1-t),0\right)$ is a geodesic of $\hbox{\tt M}$ (see
\cite{oneill}).
In the following, we shall study directly
$\N$ rather than
$\T$, since our conclusions could be easily pushed down with respect to the action of
$D$.
Consider now the holomorphic Riemannian manifold
$
\M=\left[\CI^2
\setminus ((1,i)\CI \cup (1,-i)\CI),
{du\odot dv}/({u^2+v^2})
   \right]
$.

\begin{lemma}
The geodesic equations of
both $\M$ and $\N$ are:
$
\UI^{\bullet\bullet}=
{2u}/({u^2+v^2})
\UI^{\bullet}{}^2$,
$
\VI^{\bullet\bullet}=
{2v}/({u^2+v^2})
\VI^{\bullet}{}^2
$; they are meant to be real or complex
depending on the fact that they concern $\hbox{\tt M}$ or $\hbox{\tt N}$.
\labelle{geodeq}
\end{lemma}
\begin{proposition}
All null geodesics of $\N$ are complete.
\vskip0.1truecm
\noindent
{\bf Proof:} We may deal
with the only case
$v=const:=A$.
Lemma \ref{geodeq} imply
$\UI^{\bullet\bullet}=
{2u}/({u^2+A^2})
\UI^{\bullet}{}^2$, which is solved by
$t\mapsto (C-Bt)^{-1}$ if $A=0$ and by
$t\mapsto\tan(At+B)$ if $A\not=0$, for suitable real constants $B$ and $C$.
The above functions are restrictions of meromorphic functions, hence, by definition \ref{completessa}, yield complete
geodesics.
\QUAN
\end{proposition}
We turn to nonnull geodesics of $\N$:
\begin{lemma}
The Cauchy's problem
$
\FI^{\bullet}=2A\,\Ch\varphi
\sqrt{B^2-2/A\,\Ch\varphi}$
$ \varphi(0)=\varphi_0$,
(with $B^2-2/A\,\Ch\varphi_0>0$)
has complete solutions, in the real domain, with respect to the canonical complexification, if and only if $0<AB^2\leq 2$.
\labelle{fit}
\end{lemma}
{\bf Proof:}
Set
$F(\varphi)=2A\,\Ch\varphi
\sqrt{B^2-2/A\,\Ch\varphi}$
and
$G(\varphi):=
\int_{\varphi_0}^{\varphi}
d\,\nu/F(\nu)
$, where by the integral sign
we mean the choice of the only primitive of $1/F$ vanishing at $\varphi_0$.
Rewrite the problem in the form $G(\varphi)=
\id$: this shows that $\varphi$
and $G$ are inverse elements of holomorphic functions in neighbourhoods of $\varphi_0$
and $G(\varphi_0)$.

Suppose $AB^2\geq 2$ or $AB^2<0$: then
$F$ never vanishes; since $1/F(\nu)=O(\e^{-\vert\nu\vert})$ as $\nu\to\infty$,
$G$ takes a bounded set of values, hence, by lemma \ref{inverse}, $\varphi$ is not complete.

If, instead, $0<AB^2\leq 2$, then
there exists a branch of $F$ admitting a zero on the real line, hence there exists a branch $\tilde f$ of $1/F$ whose absolute value takes all large enough values.
However $\tilde f$ can be analytically continued, by admitting complex trips, up to
$\{\varphi:\Ch\varphi\geq 2/AB^2\}$, in such a way that an even function $f$ is yielded.

Now $\vert\int_{\varphi_0}^{\varphi}f(\nu)\,d\nu\vert$ takes {\sf all} {positive} values; but $g:=\int_{\varphi_0}^{\varphi}f$ is an odd function plus a real constant
on $\{\varphi:\Ch\varphi\geq 2/AB^2\}$, hence it takes {\sf all} real values with at most the exception of its asympotical value $\sigma$.
Thus, if $(S,\varrho,\ell,H)$ is the Riemann surface of $\varphi$, then,
by lemma \ref{inverse},
$\varrho(H^{-1}(\ERRE))\cap\ERRE\supset g(\ERRE)\supset\ERRE
\setminus\{\sigma\}$.
\QUAN

\begin{definition}
The {\sf impulse function}
$\P:T\N\setminus\{{\tt null\  vectors}\}\rightarrow\ERRE$ is defined by setting $\P(\alpha,\beta,x,y)=
(\alpha^2+\beta^2)^{-1}
(2\alpha\beta+\alpha^2y/x+\beta x/y)$.
\end{definition}

\begin{theorem}
A nonnull geodesic $\gamma$ starting
from $(\alpha,\beta)$, with velocity
$(x,y)$ is complete if and only
if $0<\P(\alpha,\beta,x,y)\leq 2$.
\labelle{principal}
\end{theorem}
{\bf Proof:}
We may suppose
$\alpha\not=0$ and
$\beta\not=0$.
Moreover, we have $x\not=0$ and
$y\not=0$.
The equations in lemma \ref{geodeq}
can be integrated once to yield:
\begin{equation}
\UI^{\bullet}\VI^{\bullet}=
A(u^2+v^2),\quad u/\UI^{\bullet}+v/\VI^{\bullet}=B
,
\labelle{intprimm}
\end{equation}
where $A=xy/(\alpha^2+\beta^2)$
and $B=\alpha/x+\beta/y$; note that $AB^2=\P(\alpha,\beta,x,y)$.

Introduce now the supplementary hypothesis that $u>0$ and $v>0$:
by performing the change of coordinates $u=\e^{\omega}$,
$v=\e^{\eta}$, (\ref{intprimm})
is turned into
\begin{equation}
\oi^{\bullet}\ei^{\bullet}=
2A\, \Ch(\omega-\eta),\quad
1/\oi^{\bullet}+1/\ei^{\bullet}=B
.
\labelle{intnoeuv}
\end{equation}
We can solve with respect to
$\oi^{\bullet}$ and $\ei^{\bullet}$, getting
\begin{equation}
\cases{
\oi^{\bullet}=
2\left(B-\sqrt{B^2-2/[A\, \Ch(\oi-\ei)]}
   \right)^{-1}\cr
\ei^{\bullet}=
2\left(B+\sqrt{B^2-2/[A\, \Ch(\oi-\ei)]}
   \right)^{-1}
.
}
\labelle{intnoeuv2}
\end{equation}
Subtract and set $\varphi:=\oi-\ei$; this yields the equation in $\varphi$ studied in lemma \ref{fit},
with the appropriate initial value $\varphi(0)=\log(u/v)$; this Cauchy's problem has complete solutions
if and only if $0<\P(\alpha,\beta,x,y)\leq 2$.

Now the fact that $\varphi$ is incomplete easily implies that
so is $\gamma$.
Suppose, instead, that $\varphi$
is complete: from (\ref{intnoeuv2}),
we get
that bothm
$\ei^{\bullet}$ and $\oi^{\bullet}$ is complete; since
passing to a primitive preserves completeness,
so are $\ei$ and $\oi$:
but $u=e^{\oi}$ and
$v=e^{\ei}$: this eventually implies that
$\gamma$ is complete.

To remove the hypothesis that $u>0$ and $v>0$, consider two geodesics $\gamma$, $\delta$, starting from, say, $(\alpha,0   )$, the former with velocity
$(x,y)$ and the latter $(x,-y)$ ($y>0$).
The first order systems, like (\ref
{intprimm}), of $\gamma$ and $\delta$ differ only in the signs of constants in their first equations.
Thus, the equations of those pieces of $\gamma$
lying in $Q_1=\{u>0,v>0\}$ and of those ones of $\delta$ lying in $Q_2=\{u>0,v<0\}$ are transformed into the same system
(\ref{intnoeuv}) by performing the change of coordinates
$(u,v)=(e^{\omega}, e^{\eta})$
in $Q_1$
, resp. $(u,v)=(e^{\omega}, -e^{\eta})$ in $Q_2$; an analogous argument holds for the other octants.
It is easily seen that if a nonnull geodesic intersects one of the coordinate axes at a point $P$, it does with finite (nonnull) velocity, hence it can be analytically continued across $P$, changing octant:
thus, once obtained the (maximal) curve
$t\mapsto
(\omega(t),\eta(t))$, we can reconstruct
the original (maximal) geodesic
$t\mapsto
(u(t),v(t))$
by choosing the only smooth curve starting from
$(\alpha,\beta)$ whose graph is contained in the set
$(t,u,v\in\ERRE^3): u=\pm \e^{\omega(t)},v=\pm\e^{\eta(t)}$.
\QUAN

\label{lastpage}
\ 
\\
\tt
The author's adress:\\
CLAUDIO MENEGHINI\\
FERMO POSTA CHIASSO 1\\
CH 6830 CHIASSO - SWITZERLAND\\

\end{document}